\documentclass[aps, pre, twocolumn]{revtex4-1}
\usepackage[T1]{fontenc}
\usepackage[latin9]{inputenc}
\usepackage{amsmath}
\usepackage{graphicx}
\usepackage{stackengine}
\usepackage[caption=false]{subfig}



\begin{document}

\date{\today}
\title{Optimizing Synchronization Stability of the Kuramoto Model in Complex Networks and Power Grids}
\author{Bo Li and K. Y. Michael Wong}
\affiliation{Department of Physics, The Hong Kong University of Science and Technology, Hong Kong}

\begin{abstract}
Maintaining the stability of synchronization state is crucial for
the functioning of many natural and artificial systems. In this study,
we develop methods to optimize the synchronization stability of the Kuramoto
model by minimizing the dominant Lyapunov exponent. Using the recently proposed cut-set space approximation of the steady
states, we greatly simplify the objective function, and further derive
its gradient and Hessian with respect to natural frequencies, which
leads to an efficient algorithm with the quasi-Newton's method. The optimized
systems are demonstrated to achieve better synchronization stability
for the Kuramoto model with or without inertia in certain regimes. Hence
our method is applicable in improving the stability of power grids. It
is also viable to adjust the coupling strength of each link to improve
the stability of the system. Various operational constraints can also
be easily integrated into our scope by employing the interior point
method in convex optimization. The properties of the optimized networks
are also discussed.
\end{abstract}


\maketitle

\section{Introduction}

Synchronization occurs widely in many natural and artificial systems,
such as firefly flashes, pacemaker cells of heart, Josephson junctions, 
and power grids \cite{strogatz2004sync,pikovsky2003synchronization,Acebron2005,Kundur1994}.
In general, the synchronous states are subject to different kinds
of perturbations, and maintaining the stability of the systems against
these perturbations is crucial for the functioning of the systems
under consideration. For instance, the power grids are subject to
various disturbances and real time active controls are needed to maintain
a stable synchronization state \cite{Kundur1994}. The future power
grids will sustain larger and larger fluctuations with the introduction
of more and more renewable energies such as wind and solar power,
which raise needs to enhance the robustness and stability of existing
power networks \cite{Carrasco2006}. 

To describe these synchronization phenomena, statistical physicists
have proposed many simple but explanatory models, e.g., chaotic oscillator systems, the Kuramoto 
model, and their various generalizations \cite{Boccaletti2002,Arenas2008, Acebron2005,Rodrigues2015}. A remarkable relation between spectral aspects of network structure and synchronizability in a broad range of coupled oscillator models has been developed in the master stability function (MSF) framework \cite{Pecora1998, Arenas2008}. In particular, the second smallest eigenvalues of the graph Laplacian matrix $\lambda_2$, namely the graph algebraic connectivity, is crucial in the synchronizability of models with unbounded MSF \cite{Arenas2008}. The graph algebraic connectivity is an interesting measure of network connectivity \cite{DEABREU2007, chungspectral}, whose role in dynamical stability can be exemplified in consensus dynamics or diffusion on networks $\dot{x}_i = -\sum_j L_{ij} x_j$, where $\lambda_2$ determines the rate of convergence of the slowest mode \cite{Arenas2008}. The graph algebraic connectivity is solely determined by the network topology. However, in many networks such as the power grid and transportation networks, stable behavior also depends on attributes other than topology.

In this study, we focus on the stability of the Kuramoto model on general networks. Due to the heterogeneity of power supply and demand, the stability of the frequency synchronization state of this nonlinear dynamical model is no longer determined by the graph algebraic connectivity or network structure itself, but is replaced by an algebraic connectivity that has an intricate dependence on the system steady state \cite{ARAPOSTHATIS1981}. The optimization of synchronization stability should take into account both the graph connectivity and the dynamical parameters.

Enhancing the synchronization stability in these settings has been stressed in a few recent studies \cite{Mallada2011, Motter2013}, where the effects of network structures
or power grid parameters, e.g., the damping coefficients and power
injections, on the system stability were explored. However, a practical consideration in implementing real-time flow control of the networks is the efficiency in calculating the gradient of the objective function in the space of variables, as was done in the cases of power scheduling and line impedance modification in power grids. Conventionally this requires us to solve the nonlinear flow equations in each update step, seriously slowing down the process. In this paper, we introduce the cut-set space approximation~\cite{Dorfler2013}, enabling us to express the objective function in terms of the graph algebraic connectivity, thereby saving the need for the stepwise solution of the nonlinear flow equation and greatly simplifying the calculation of power flow and the evaluation of the gradients of the objective function.

\section{The Model}

\subsection{First-order Kuramoto model}

We focus on the non-uniform first-order Kuramoto model on a connected
network in the form of 
\begin{equation}
\dot{\theta}_{i}=\omega_{i}+\sum_{j}K_{ij}\sin(\theta_{j}-\theta_{i}),\label{eq:1st_order_model}
\end{equation}
where $\theta_{i}$ denotes the phase angle of node $i$, $\omega_{i}$
the natural frequency and $K_{ij}(=K_{ji})$ the coupling strength between
node $i$ and node $j$. Without loss of generality, we assume $\sum_{i}\omega_{i}=0$.
The steady state is given by 
\begin{equation}
0=\omega_{i}+\sum_{j}K_{ij}\sin(\theta_{j}^{*}-\theta_{i}^{*}).\label{eq:steady_state}
\end{equation}

In the leading order, the small deviation from the steady state $\delta\theta_{i}=\theta_{i}-\theta_{i}^{*}$
follows \cite{ARAPOSTHATIS1981}
\begin{align*}
\delta\dot{\theta}_{i} & \approx\sum_{j}K_{ij}\cos(\theta_{j}^{*}-\theta_{i}^{*})(\delta\theta_{j}-\delta\theta_{i})\\
 & =-\sum_{j}L(\theta^{*})_{ij}\delta\theta_{j},
\end{align*}
where $L(\theta^{*})_{ij}:=\delta_{ij}\sum_{l}K_{il}\cos(\theta_{l}^{*}-\theta_{i}^{*})-K_{ij}\cos(\theta_{j}^{*}-\theta_{i}^{*})$
is a state-dependent Laplacian matrix with edge weight $W(\theta^{*})_{ij}=K_{ij}\cos(\theta_{j}^{*}-\theta_{i}^{*})$.
Note that this Laplacian matrix depends on the steady state of the system, in contrast with the state-independent Laplacian, which we denote as $L[K]_{ij}:=\delta_{ij}\sum_l K_{il} - K_{ij}$. The Jacobian matrix is $J=-L(\theta^{*})$, which has a null-space
of dimension one, corresponding to the rotational symmetry of the
model. If $|\theta_{j}^{*}-\theta_{i}^{*}|<\pi/2$ holds for every
edge $(i,j)$, then all the edge weights $W_{ij}$ are positive and the lowest eigenvalue is 0, corresponding to the mode of uniform displacement. All the other eigenvalues of $L(\theta^{*})$ are positive, making
the dynamical system locally exponentially stable. In this case, the
slowest mode corresponds to the second lowest eigenvalue of $L(\theta^{*})$,
that is, the negative of the largest Lyapunov exponent excluding the
null exponent of $J$. We denote it as $\lambda_{2}(L(\theta^{*}))$
and call it the \textit{state algebraic connectivity} to distinguish it from
the usual \textit{graph algebraic connectivity} $\lambda_{2}(L[K])$. To improve
the stability, our objective is to maximize $\lambda_{2}(L(\theta^{*}))$
as in Ref. \cite{Mallada2011}.

\subsection{Second-order Kuramoto model}

The second-order Kuramoto model is gaining attention due to its resemblance
to the swing equation of power grids neglecting the transmission losses
\cite{ARAPOSTHATIS1981}

\begin{equation}
M_{i}\ddot{\theta}_{i}+D_{i}\dot{\theta}_{i}=P_{i}+\sum_{j}\frac{|V_{i}||V_{j}|}{X_{ij}}\sin(\theta_{j}-\theta_{i}),
\end{equation}
where $M_{i}$ and $D_{i}$ are the inertia and damping coefficient
of node $i$ respectively, $P_{i}$ and $|V_{i}|$ are the mechanical power and
voltage magnitude of node $i$, and $X_{ij}$ is the line reactance
of edge $(i,j)$. The connection to the Kuramoto model is obvious if $P_{i}$
is identified as the natural frequency $\omega_{i}$ and $|V_{i}||V_{j}|/X_{ij}$
is identified as coupling $K_{ij}$. For simplicity, we consider uniform
inertia and damping coefficient $M_{i}=M$ and $D_{i}=D$ and focus
on the following model 
\begin{equation}
M\ddot{\theta}_{i}+D\dot{\theta}_{i}=\omega_{i}+\sum_{j}K_{ij}\sin(\theta_{j}-\theta_{i}).\label{eq:2nd_order_model}
\end{equation}

The steady state $(\dot{\theta}^{*}=0,\theta^{*})$ is again given
by Eq. (\ref{eq:steady_state}), with the Jacobian matrix evaluated at
this point as \cite{Mallada2011,Motter2013}

\[
J(\dot{\theta}^{*}=0,\theta^{*})=\begin{bmatrix}-\frac{D}{M}I & -\frac{1}{M}L(\theta^{*})\\
I & 0
\end{bmatrix}.
\]

As derived in Ref. \cite{Motter2013}, $J(\dot{\theta}^{*},\theta^{*})$
can be diagonalized by the eigenvectors of $L(\theta^{*})$, with
corresponding eigenvalues
\begin{equation}
\mu_{j\pm}(\lambda_{j},D,M)=-\frac{D}{2M}\pm\frac{1}{2}\sqrt{\left(\frac{D}{M}\right)^{2}-\frac{4}{M}\lambda_{j}(L(\theta^{*}))}.
\end{equation}
The maximal nontrivial eigenvalue is $\mu_{2+}=-\frac{D}{2M}+\frac{1}{2}\sqrt{\left(\frac{D}{M}\right)^{2}-\frac{4}{M}\lambda_{2}(L(\theta^{*}))}.$
When $\lambda_{2}(L(\theta^{*}))<D^{2}/4M$, improving $\lambda_{2}(L(\theta^{*}))$
will always lead to the increment of $\mu_{2+}$. In this regime, optimizing $\lambda_{2}(L(\theta^{*}))$
is also applicable to stabilizing the uniform second-order Kuramoto
model, therefore it can be applied in the stabilization of power grids. This regime can correspond to large damping, small
inertia or close to bifurcation.

\section{Method}

\subsection{Variation of state algebraic connectivity}

Viewing $\omega_{i}$ and $K_{ij}$ as control variables, we aim at
maximizing $\lambda_{2}(L(\theta^{*}))$ in order to improve the stability
of both Eqs. (\ref{eq:1st_order_model}) and (\ref{eq:2nd_order_model}).
We first derive the variation of state algebraic connectivity due
to change of natural frequency. We assume that the state algebraic
connectivity is non-degenerate throughout optimization, which usually
holds when the corresponding graph algebraic connectivity $\lambda_{2}(L[K])$
is non-degenerate.

There is no explicit expression of $\lambda_2(L(\theta^*))$. Nevertheless, it is possible to derive its derivatives using the perturbation theory, as commonly practiced in quantum mechanics. In the case that $\lambda_2(L_(\theta^*))$ is non-degenerate, the variation of $\lambda_{2}(L(\theta^{*}))$
is given by \cite{landau1977quantum}
\begin{align}
\delta\lambda_{2}(L(\theta^{*})) & =\langle v_{2}(\theta^{*})|\delta L(\theta^{*})|v_{2}(\theta^{*})\rangle\nonumber \\
 & =v_{2}(\theta^{*})^{T}\,\delta L(\theta^{*})\,v_{2}(\theta^{*}),
\end{align}
where $v_{2}(\theta^{*})$ is the normalized eigenvector of $L(\theta^{*})$
corresponding to $\lambda_{2}(L(\theta^{*}))$. Since $L(\theta^{*})$
is a Laplacian matrix with edge weight $W(\theta^{*})_{ij}=K_{ij}\cos(\theta_{j}^{*}-\theta_{i}^{*})$,
one has $\delta L(\theta^{*})_{ij}=\delta_{ij}\sum_{l}\delta W(\theta^{*}){}_{il}-\delta W(\theta^{*}){}_{ij}$
and
\begin{equation}
\delta\lambda_{2}(L(\theta^{*}))=\sum_{(i,j)}\delta W(\theta^{*}){}_{ij}[v_{2}(\theta^{*})_{i}-v_{2}(\theta^{*})_{j}]^{2}.
\end{equation}
So the gradient of the state algebraic connectivity with respect to $\omega$
is
\begin{equation}
[\nabla_{\omega}\lambda_{2}(L(\theta^{*}))]_{k}=\sum_{(i,j)}\frac{\delta W(\theta^{*})_{ij}}{\delta\omega_{k}}[v_{2}(\theta^{*})_{i}-v_{2}(\theta^{*})_{j}]^{2}.\label{eq:gradient_omega}
\end{equation}

The computational complexity comes from the implicit dependence between shift
of steady state $\delta\theta^{*}$ and change of natural frequency
$\delta\omega$. In Ref. \cite{Mallada2011}, $\delta\theta^{*}/\delta\omega$
is proved to be related to the pseudo-inverse of $L(\theta^{*})$.
These expressions lead to a gradient ascent method to maximize $\lambda_{2}(L(\theta^{*}))$
by scheduling $\omega$. However, this method requires solving the
steady state equation Eq. (\ref{eq:steady_state}) and computing the
pseudo-inverse of $L(\theta^{*})$ in every iteration, both of which
are time consuming. In addition, convergence to the optimal solution
can be very slow for gradient ascent update. In this paper, we propose
to use the cut-set space approximation to simplify the problems as
follows.

\subsection{Cut-set space approximation of network flows}

The natural frequency $\omega_{i}$ can be viewed as supply or demand of node $i$ in a supply
network as $P_{i}$ in the power grid, and $K_{ij}\sin(\theta_{j}-\theta_{i})$
is the resource or power transported from node $j$ to node $i$. The
steady-state Eq. (\ref{eq:steady_state}) implies the flow
conservation on each node. 

Solving the nonlinear steady-state equation can be computationally costly.
Recently, it has been shown that the cut-set space approximation of
power flows can be rather accurate in many regimes \cite{Dorfler2013,DorflerIEEE2013}. For completeness, the main steps are outlined as follows.
We first formally rewrite the anti-symmetric quantity $\sin(\theta_{j}^{*}-\theta_{i}^{*})$ as $\beta_{ij}(=-\beta_{ji})$, which we try to decompose into the sum of two parts $\beta_{ij} = \beta_{ij}^{\text{cut}} + \beta_{ij}^{\text{cycle}}$. The first part $\beta_{ij}^{\text{cut}}$ is expressed by the potential difference $\beta_{ij}^{\text{cut}}=\phi_{j}-\phi_{i}$, where $\phi_i$ is an unknown potential function to be solved self-consistently. The second part $\beta_{ij}^{\text{cycle}}$ satisfies the circular flow relation $\sum_{j\in\partial i}K_{ij}\beta_{ij}^{\text{cycle}}=0\, \forall i$. In the language of graph theory, $\beta^{\text{cut}}$ and $\beta^{\text{cycle}}$ are said to live in the cut-set space and cycle space respectively \cite{Biggs1997, Dorfler2013}. Substituting $\beta_{ij} = \phi_j - \phi_i + \beta_{ij}^{\text{cycle}}$ into Eq. (\ref{eq:steady_state}), we have
\begin{align}
0 & = \omega_i + \sum_{j\in\partial i} K_{ij} (\phi_j - \phi_i + \beta_{ij}^{\text{cycle}}) \nonumber \\
  & = \omega_i - \sum_{j\in\partial i} L[K]_{ij} \phi_j, \label{eq:beta_cut_eq_of_state}
\end{align}
where $L[K]$ is the graph Laplacian matrix, which depends only on the network topology and edge weights. By taking the pseudo-inverse of $L[K]$, denoted as $L[K]^{\dagger}$, the potential $\phi$ is obtained by $\phi = L[K]^{\dagger}\omega$, and subsequently, $\beta_{ij}^{\text{cut}}=\phi_{j}-\phi_{i}=\sum_{l}(L[K]_{jl}^{\dagger}-L[K]_{il}^{\dagger})\omega_{l}$. It turns out that $\phi$ coincides
with the DC approximation of AC power flow in power engineering $\theta^{\text{DC}}$
\cite{Dorfler2013}. To simplify
the calculation, it is proposed to approximate $\beta$ by its cut-set
space component $\beta^{\text{cut}}$, i.e., $\sin(\theta_{j}^{*}-\theta_{i}^{*})\approx\phi_{j}-\phi_{i}=\sum_{l}(L[K]_{jl}^{\dagger}-L[K]_{il}^{\dagger})\omega_{l}$.

\begin{figure}
\hspace{-6mm}
\subfloat{
\topinset{(a)}{\includegraphics[scale=0.21]{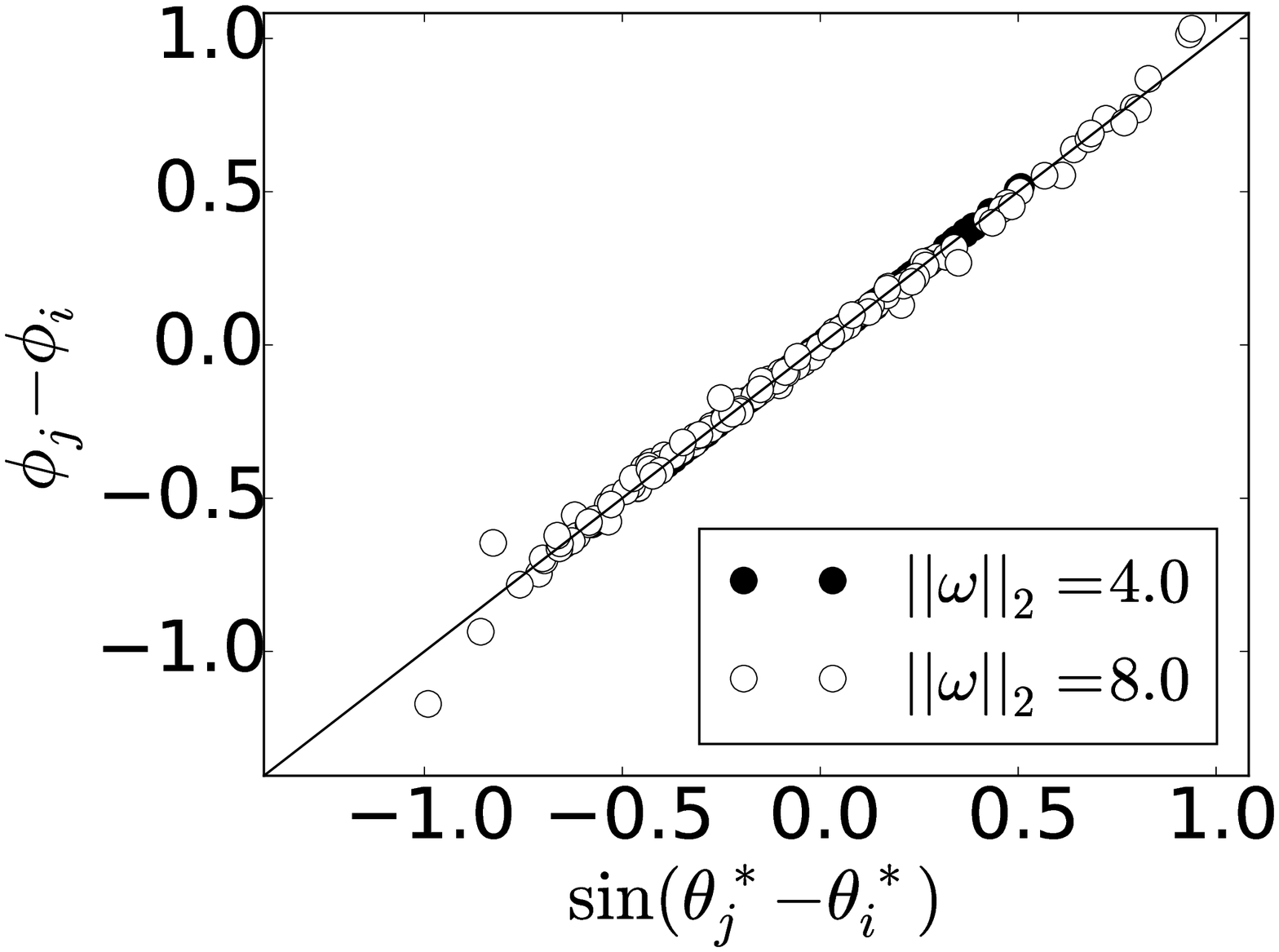}}{1mm}{-19mm}
\begin{picture}(0,0)
\put(-99,58){\includegraphics[scale=0.0615]{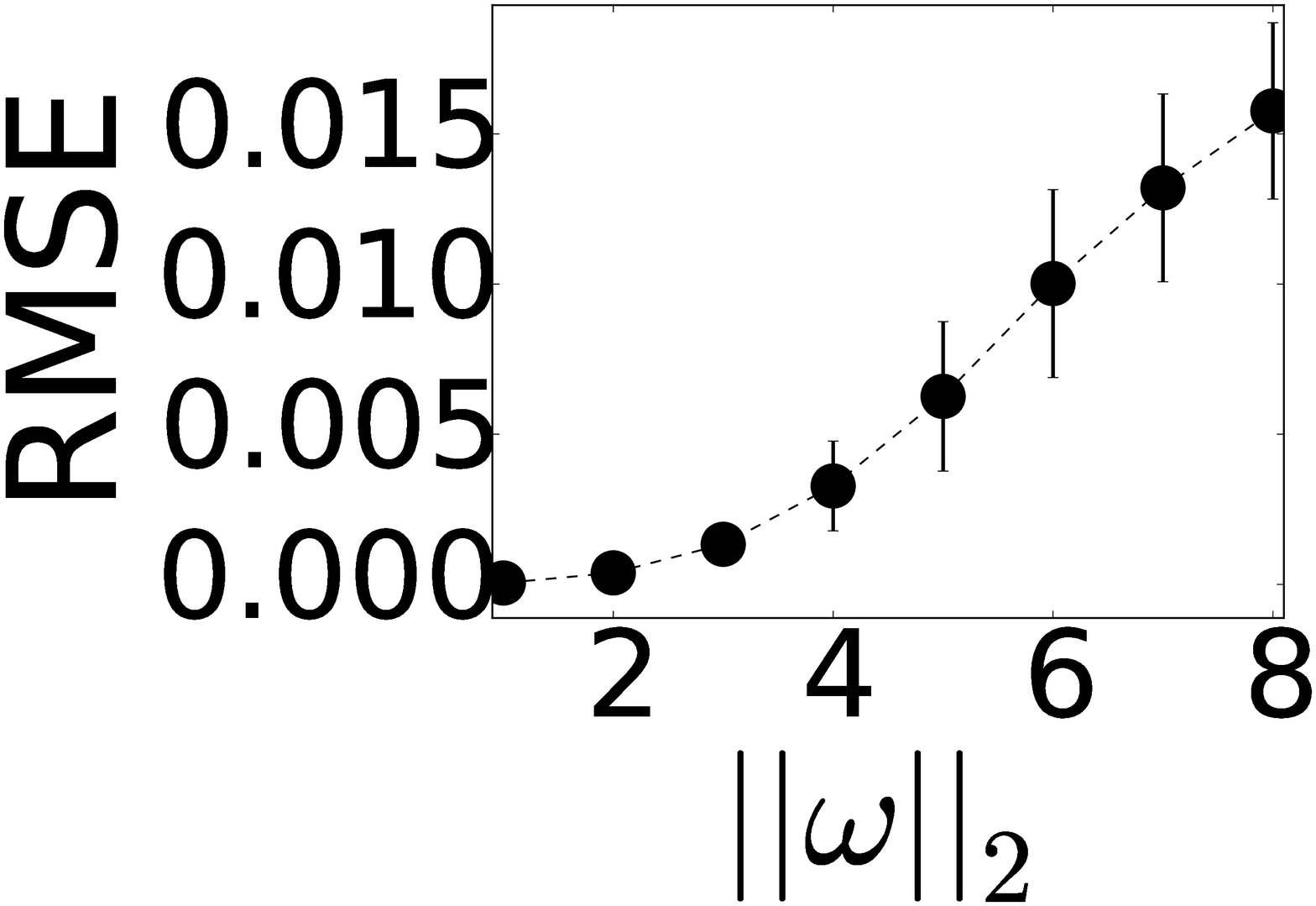}}
\end{picture}
}
\hspace{-4mm}
\subfloat{
\topinset{(b)}{\includegraphics[scale=0.21]{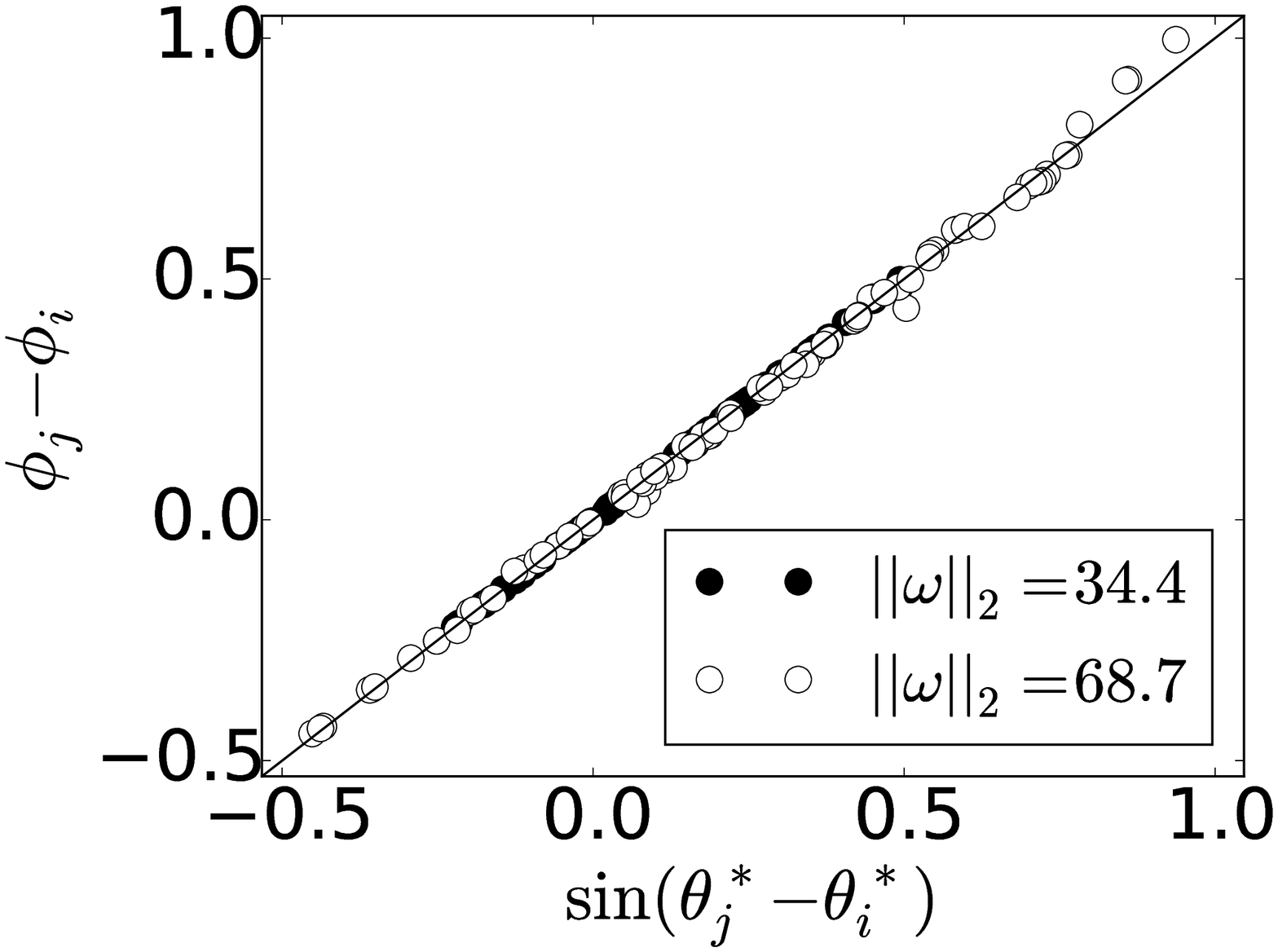}}{1mm}{-19mm}
\begin{picture}(0,0)
\put(-99,58){\includegraphics[scale=0.0615]{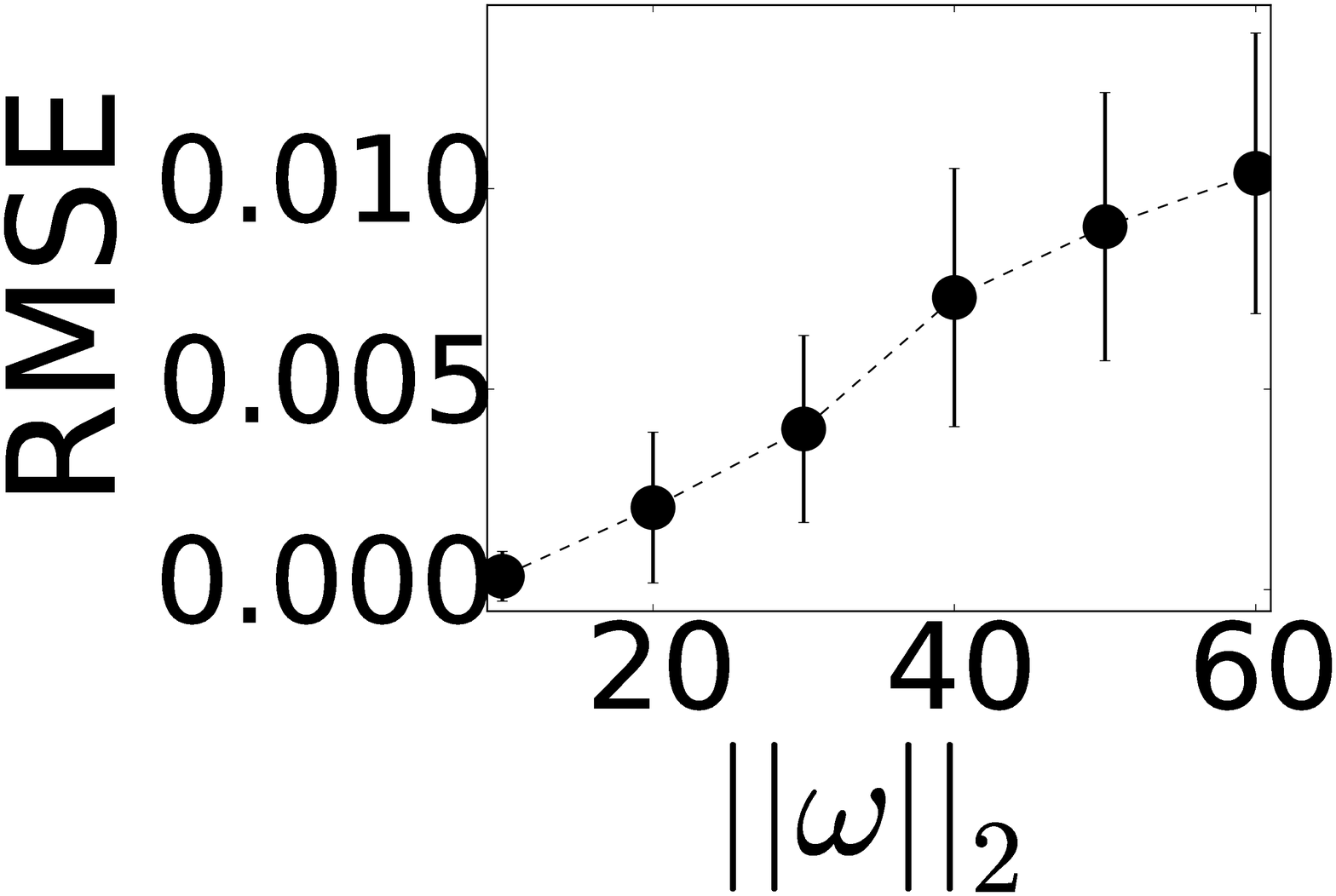}}
\end{picture}
}

\caption{$\phi_{j}-\phi_{i}$ vs. $\sin(\theta_{j}^{^{*}}-\theta_{i}^{*})$.
(a) Erd{\"o}s-R{\'e}nyi graph of 50 nodes (ER50), where $\omega$ is drawn
from a Gaussian distribution and $K_{ij}=1$. Inset: root-mean-square error (RMSE) of estimator $\phi_j - \phi_i$ for $\sin(\theta^*_j - \theta^*_i)$ among all the edges. Each data point is averaged over 100 samples.  (b) IEEE reliability test
system 96 (RTS96) \cite{Grigg1999}, where $\omega$ is modified from the power injection
data in the test system and $K_{ij}$ is defined to be the inverse
of line reactance of edge $(i,j)$. Inset: RMSE of estimator $\phi_j - \phi_i$ for $\sin(\theta^*_j - \theta^*_i)$ among all the edges. Each data point is averaged over 100 samples. \label{fig:cut_set_approx} }
\end{figure}

Such an approximation is exact in some specific systems, such as acyclic graphs and systems with cut-set inducing frequencies, while it has also been tested numerically in many generic networks that the approximation is surprisingly accurate \cite{Dorfler2013,DorflerIEEE2013}. We demonstrate
two examples in Fig.~\ref{fig:cut_set_approx}. To quantify the stress of the system, the $L^2$ norm (or the Euclidean norm) of the natural frequency is used, i.e., $\|\omega\|_2:=\sqrt{\sum_i \omega_i^2}$. It is shown that the potential difference $\phi_{j}-\phi_{i}$ approximates $\sin(\theta_{j}^{^{*}}-\theta_{i}^{*})$ quite well even in the stress cases with large $\|\omega\|_{2}$.

\subsection{Optimization by tuning natural frequencies}

With the cut-set-space approximation, the edge weight of the state-dependent
Laplacian matrix $L(\theta^{*})$ can be approximated as 
\begin{align*}
W(\theta^{*})_{ij} & =K_{ij}\cos(\theta_{j}^{*}-\theta_{i}^{*})=K_{ij}\sqrt{1-\sin^{2}(\theta_{j}^{*}-\theta_{i}^{*})}\\
 & \approx\tilde{W}(\phi)_{ij}=K_{ij}\sqrt{1-(\phi_{j}-\phi_{i})^{2}}\\
 & \equiv K_{ij}\sqrt{1-\sum_{kl}\omega_{k}A_{kl}^{(ij)}\omega_{l}}.
\end{align*}
where $A^{(ij)}$ is defined to be a matrix with entry $A_{kl}^{(ij)}=(L[K]_{jk}^{\dagger}-L[K]_{ik}^{\dagger})(L[K]_{jl}^{\dagger}-L[K]_{il}^{\dagger})$
and we have made use of the fact that $\phi=L[K]^{\dagger}\omega$.
Provided that $L[K]^{\dagger}$ is calculated and recorded, every time we calculate $W(\theta^*)$ we only
need to solve for $\phi$ by simple matrix multiplication instead of
solving the nonlinear steady-state equation Eq. (\ref{eq:steady_state}).
Now we work on the state algebraic connectivity $\lambda_{2}(\tilde{L}(\phi))$, 
which corresponds to the state-dependent Laplacian matrix with edge
weight $\tilde{W}(\phi)_{ij}=K_{ij}\sqrt{1-\sum_{kl}\omega_{k}A_{kl}^{(ij)}\omega_{l}}$.
We assume in the following discussion that $|\phi_{j}-\phi_{i}|<1$
always holds such that $\tilde{W}(\phi)_{ij}$ is real for every edge
$(i,j)$. This assumption can fail when the system is so stressed that $|\theta_{j}^{*}-\theta_{i}^{*}|$ is close to $\pi/2$ along
some edges, in which case a preprocess to destress the system before
optimization is needed.

The gradient in Eq. (\ref{eq:gradient_omega}) can be estimated by
$\nabla_{\omega}\lambda_{2}(\tilde{L}(\phi))$ 
\begin{equation}
[\nabla_{\omega}\lambda_{2}(\tilde{L}(\phi))]_{k}=\sum_{(i,j)}K_{ij}\frac{-\sum_{l}A_{kl}^{(ij)}\omega_{l}}{\sqrt{1-\omega^{T}A^{(ij)}\omega}}[v_{2}(\phi)_{i}-v_{2}(\phi)_{j}]^{2},
\end{equation}
where $v_{2}\text{(\ensuremath{\phi})}$ is the normalized eigenvector corresponding
to $\lambda_{2}(\tilde{L}(\phi))$.

Similarly, the Hessian of the state algebraic connectivity is estimated
by 
\begin{align*}
 & H_{kl}=\frac{\partial^{2}\lambda_{2}(\tilde{L}(\phi))}{\partial\omega_{k}\partial\omega_{l}}=\sum_{(i,j)}\frac{\partial^{2}\tilde{W}(\phi)_{ij}}{\partial\omega_{k}\partial\omega_{l}}[v_{2}(\phi)_{i}-v_{2}(\phi)_{j}]^{2}\\
 & +\sum_{(i,j)}2\frac{\partial\tilde{W(\phi)_{ij}}}{\partial\omega_{k}}[v_{2}(\phi)_{i}-v_{2}(\phi)_{j}]\bigg[\frac{\partial v_{2}(\phi)_{i}}{\partial\omega_{l}}-\frac{\partial v_{2}(\phi)_{j}}{\partial\omega_{l}}\bigg],
\end{align*}
where $\partial v_{2}(\phi)/\partial\omega$ can also be obtained
from the non-degenerate perturbation theory, which is computationally
costly. We found in all our numerical experiments that truncating
the second term of the Hessian can still lead to efficient optimization
but simplify the calculation significantly. Hence, in the following
we use the approximated Hessian $H_{kl}\approx\sum_{(i,j)}\partial^{2}\tilde{W}(\phi)_{ij}/\partial\omega_{k}\partial\omega_{l}[v_{2}(\phi)_{i}-v_{2}(\phi)_{j}]^{2}$
for optimization.

Obtaining the gradient and Hessian, we can define the update direction
of gradient ascent and quasi-Newton method to maximize $\lambda_{2}(\tilde{L}(\phi))$,

\begin{eqnarray*}
\Delta\omega^{\text{gradient}} & = & \nabla_{\omega}\lambda_{2}(\tilde{L}(\phi)),\\
\Delta\omega^{\text{Newton}} & = & H^{-1}\nabla_{\omega}\lambda_{2}(\tilde{L}(\phi)).
\end{eqnarray*}
The natural frequency is updated by $\omega\leftarrow\omega+s\Delta\omega^{\text{gradient}}$
or $\omega\leftarrow\omega+s\Delta\omega^{\text{Newton}}$ with the
step size $s$ determined by back tracking line search \cite{boyd2004convex},
after which $\omega$ is enforced to be zero-sum by $\omega_{i}\leftarrow\omega_{i}-1/N\sum_{j}\omega_{j}$
so that it admits a steady state.

\begin{figure}
\topinset{(a)}{\includegraphics[scale=0.22]{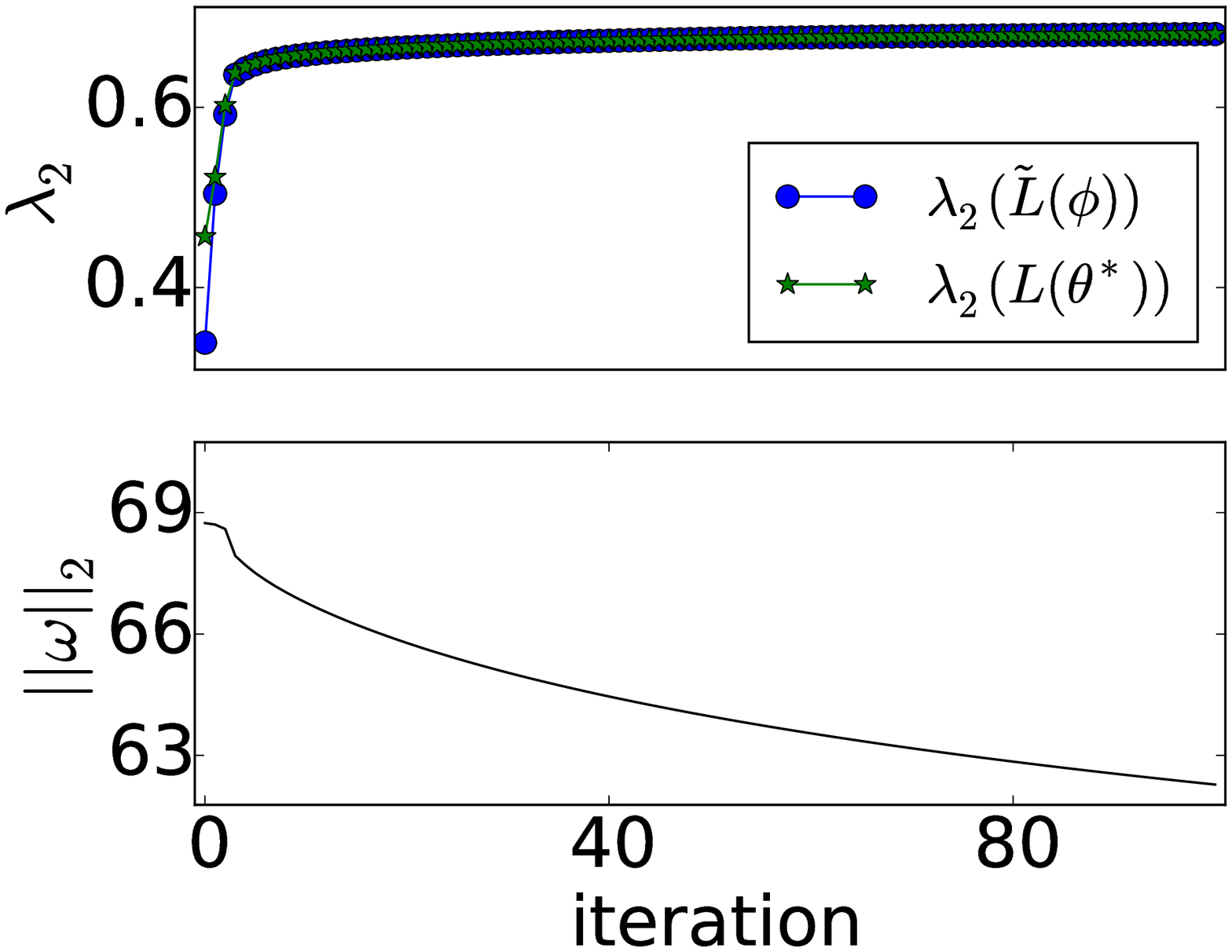}}{-1mm}{-18mm}
\topinset{(b)}{\includegraphics[scale=0.22]{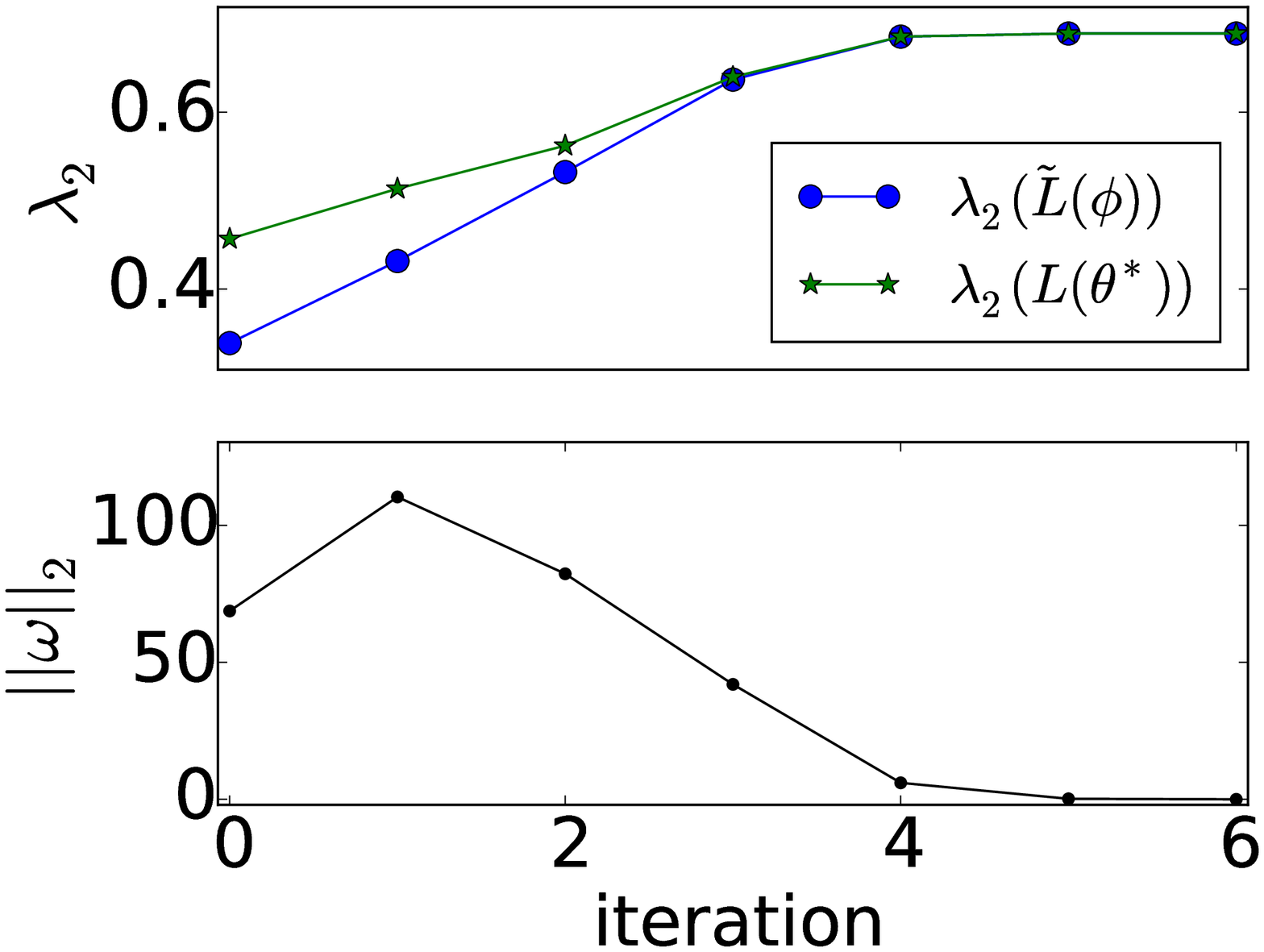}}{-1mm}{-18mm}

\caption{(Color online) $\lambda_{2}$ and $\|\omega\|_{2}$ through optimization for RTS96
power network. The initial natural frequency is modified from the
power injection data in the test case. (a) Gradient ascent update.
Both $\lambda_{2}(\tilde{L}(\phi))$ and $\lambda_{2}(L(\theta^{*}))$
increase gently in the later stage, and the natural frequency $\omega$
is approaching the optimal state $\omega=0$ very slowly due
to the flat landscape. (b) Quasi-Newton update. Both $\lambda_{2}(\tilde{L}(\phi))$
and $\lambda_{2}(L(\theta^{*}))$ approach the optimum $\lambda_{2}(L(\theta^{*}=0))=0.6889$
after six iterations. \label{fig:RTS96_unconstrained}}

\end{figure}

In general, $\lambda_{2}(\tilde{L}(\phi))$ is an increasing function
with $\tilde{W}(\phi)_{ij}$, which favors small phase angle difference
across each edge. Without imposing any constraint, the optimal solution
should take place at $\omega=0$, in which case the optimum $\lambda_{2}(L(\theta^{*}=0))$
coincides with the graph algebraic connectivity. In Fig.~\ref{fig:RTS96_unconstrained}
we show the the optimization process for the RTS96 power network with
gradient ascent update and quasi-Newton update. It is observed in
this case that (i) $\lambda_{2}(\tilde{L}(\phi))$ is close to the
exact state algebraic connectivity $\lambda_{2}(L(\theta^{*}))$ at the same $\omega$ (obtained by solving the steady state equation Eq. (\ref{eq:steady_state}) with $\omega$ given at that iteration);
(ii) the Newton's method is much more efficient than the gradient
ascent, approaching the optimum within only a few steps, despite the extra efforts for computing the Hessian $H$ and solving the linear
equation $H\Delta\omega^{\text{Newton}}=\nabla_{\omega}\lambda_{2}(\tilde{L}(\phi))$
to obtain $\Delta\omega^{\text{Newton}}$. By taking the advantages of the cut-set space approximation and the Newton's method, our approach here provides a much more efficient algorithm compared to the previous study that relied on the full calculation of the nonlinear steady state and the gradient ascent update \cite{Mallada2011}.

\subsection{Optimization by tuning for coupling strengths}

Instead of optimizing the natural frequencies, one can also tune the
coupling strengths of edges to improve the stability. In power grids,
this corresponds to the change of line reactance of each edge, which
may be implemented by tuning the transmission lines or using FACTS
devices \cite{Gotham1998}. Similarly, we can also derive the gradient
and Hessian of $\lambda(\tilde{L}(\phi))$ with respect to the coupling
strength
\begin{align}
[\nabla_{K} & \lambda_{2}(\tilde{L}(\phi))]_{(k,l)}=\sum_{(i,j)}\frac{\delta\tilde{W}(\phi)_{ij}}{\delta K_{kl}}[v_{2}(\phi)_{i}-v_{2}(\phi)_{j}]^{2}\nonumber \\
 & =\sum_{(i,j)}\bigg\{    \delta_{(i,j),(k,l)}\sqrt{1-\omega^{T}A^{(ij)}\omega}+\nonumber \\
 & \,\,\,\,\,\,\,\, \frac{1}{2}K_{ij}\frac{-\sum_{mn}\omega_{m}\frac{\partial A_{mn}^{(ij)}}{\partial K_{kl}}\omega_{n}}{\sqrt{1-\omega^{T}A^{(ij)}\omega}}\bigg\}    [v_{2}(\phi)_{i}-v_{2}(\phi)_{j}]^{2},\label{eq:gradient_wrt_K}
\end{align}
where the evaluation of $\partial A^{(ij)}/\partial K_{kl}$ relies
on the computation of $\partial L[K]^{\dagger}/\partial K_{kl}$ which
is attainable as long as the rank of $L[K]$ remains unchanged \cite{Colub1973}.
The gradient ascent update is simply given by $K\leftarrow K+s\nabla_{K}\lambda_{2}(\tilde{L}(\phi))$.
The Hessian matrix and update of Newton's method can also be obtained
straightforwardly, although the expression is extremely tedious. The
update of coupling strength renders the modification of $L[K]$ and
recalculation of $L[K]^{\dagger}$, making it much more time consuming
than the update of natural frequencies.

Although we have been dealing with the oscillatory system with sinusoidal coupling, 
we remark that the general framework developed here can also be applicable to systems 
with other coupling functions, and even other eigenvalue optimization problems, especially when 
nonlinearity comes into play and the usual semidefinite programming is not directly applicable \cite{boyd2004convex}. 

\section{Results}

\subsection{Behavior at optimal natural frequencies}
\begin{figure}
\hspace{-1mm}
\topinset{(a)}{\includegraphics[scale=0.21]{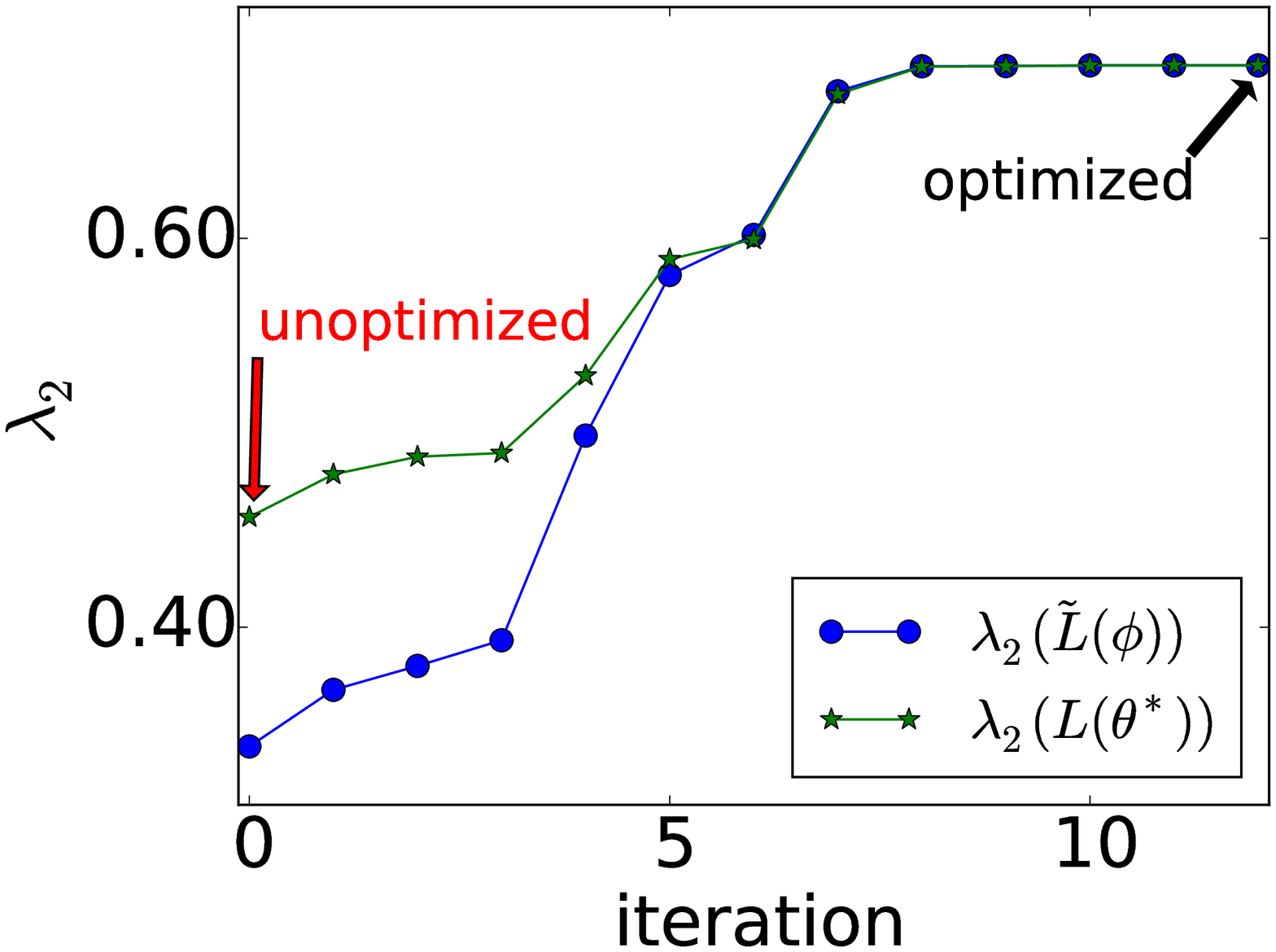}}{-1mm}{-19mm}
\topinset{(b)}{\includegraphics[scale=0.22]{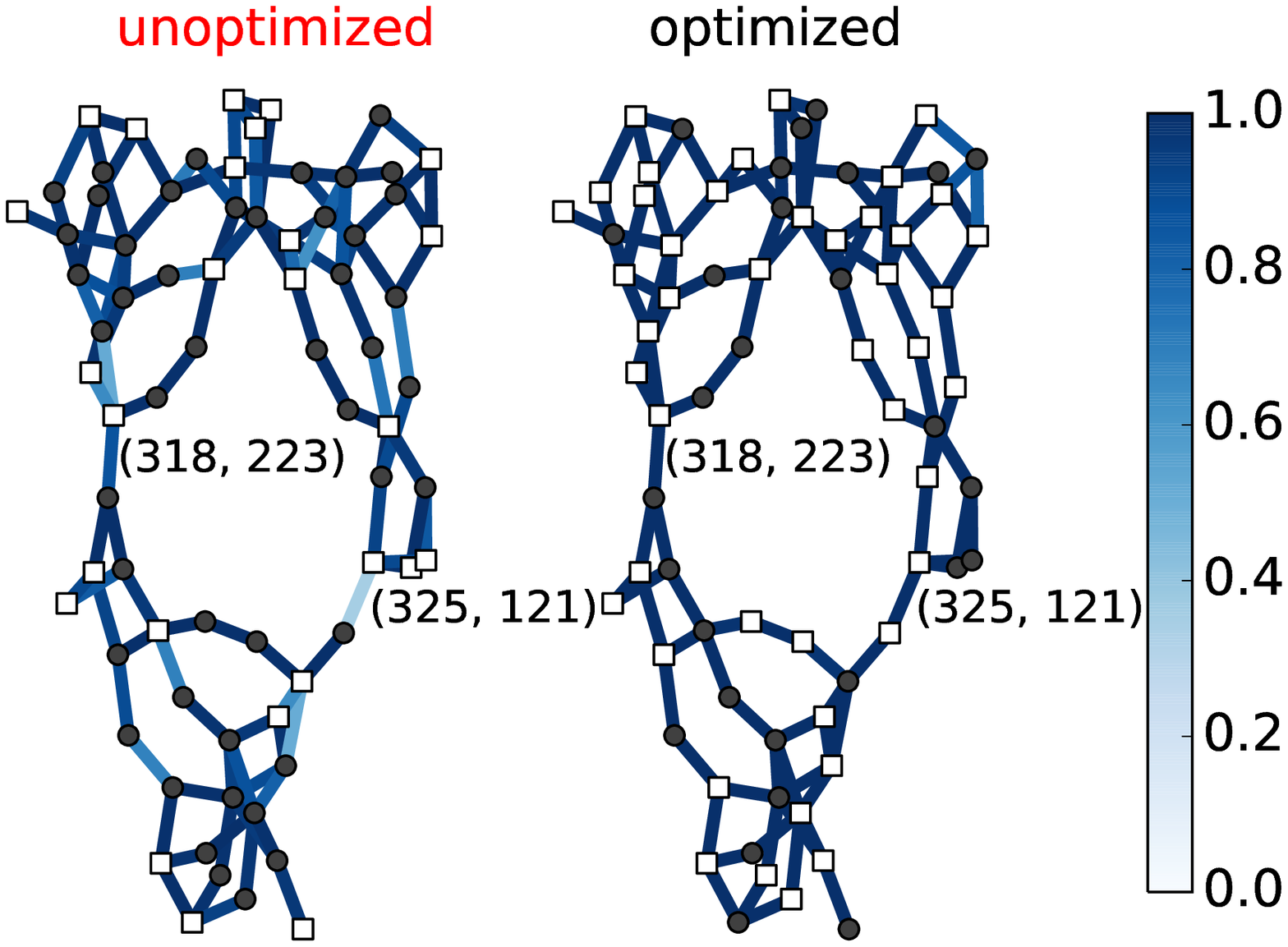}}{-1mm}{-19mm}\\
\hspace{-1mm}
\topinset{(c)}{\includegraphics[scale=0.21]{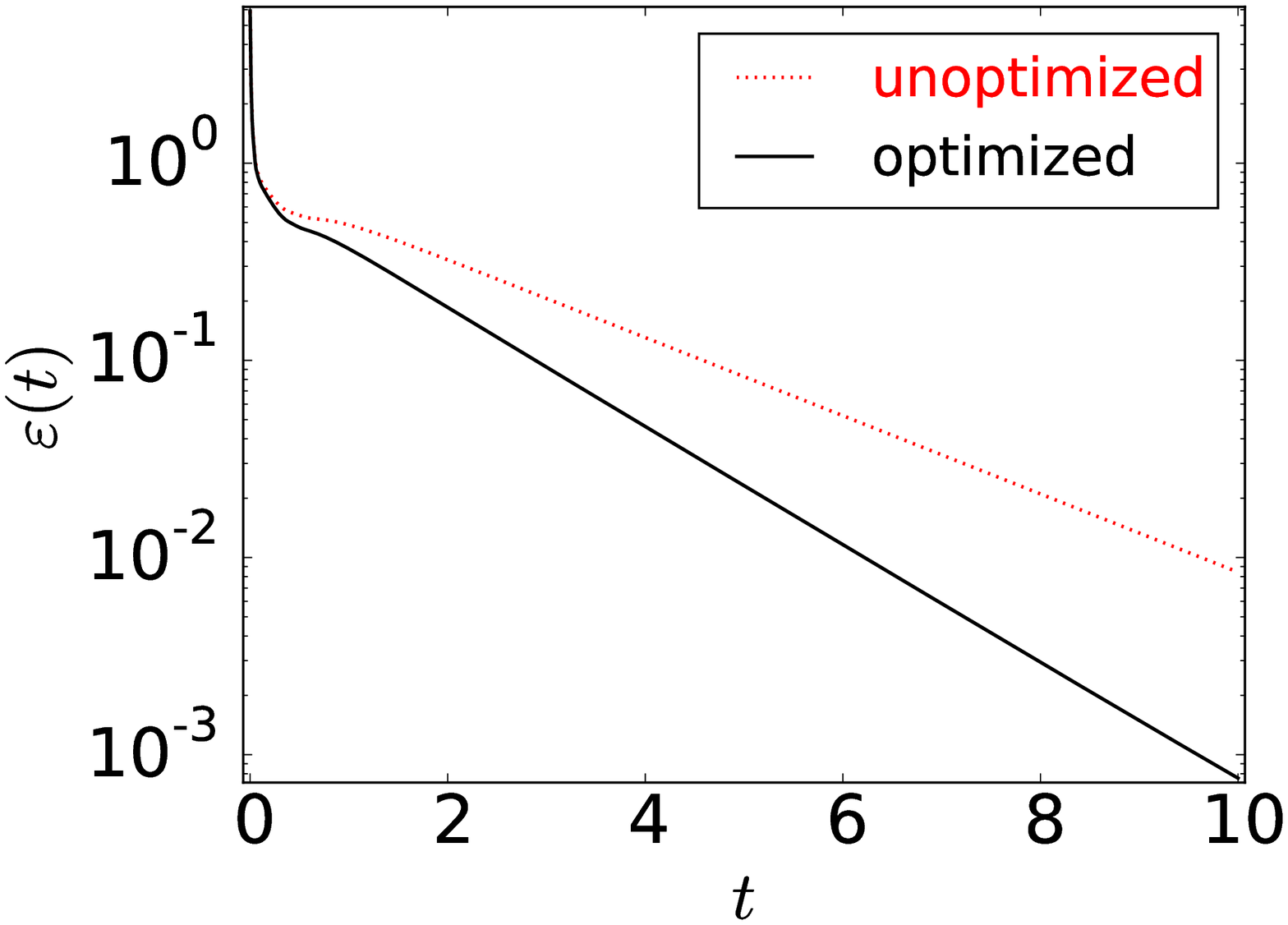}}{-1mm}{-19mm}\hspace{-1mm}
\topinset{(d)}{\includegraphics[scale=0.21]{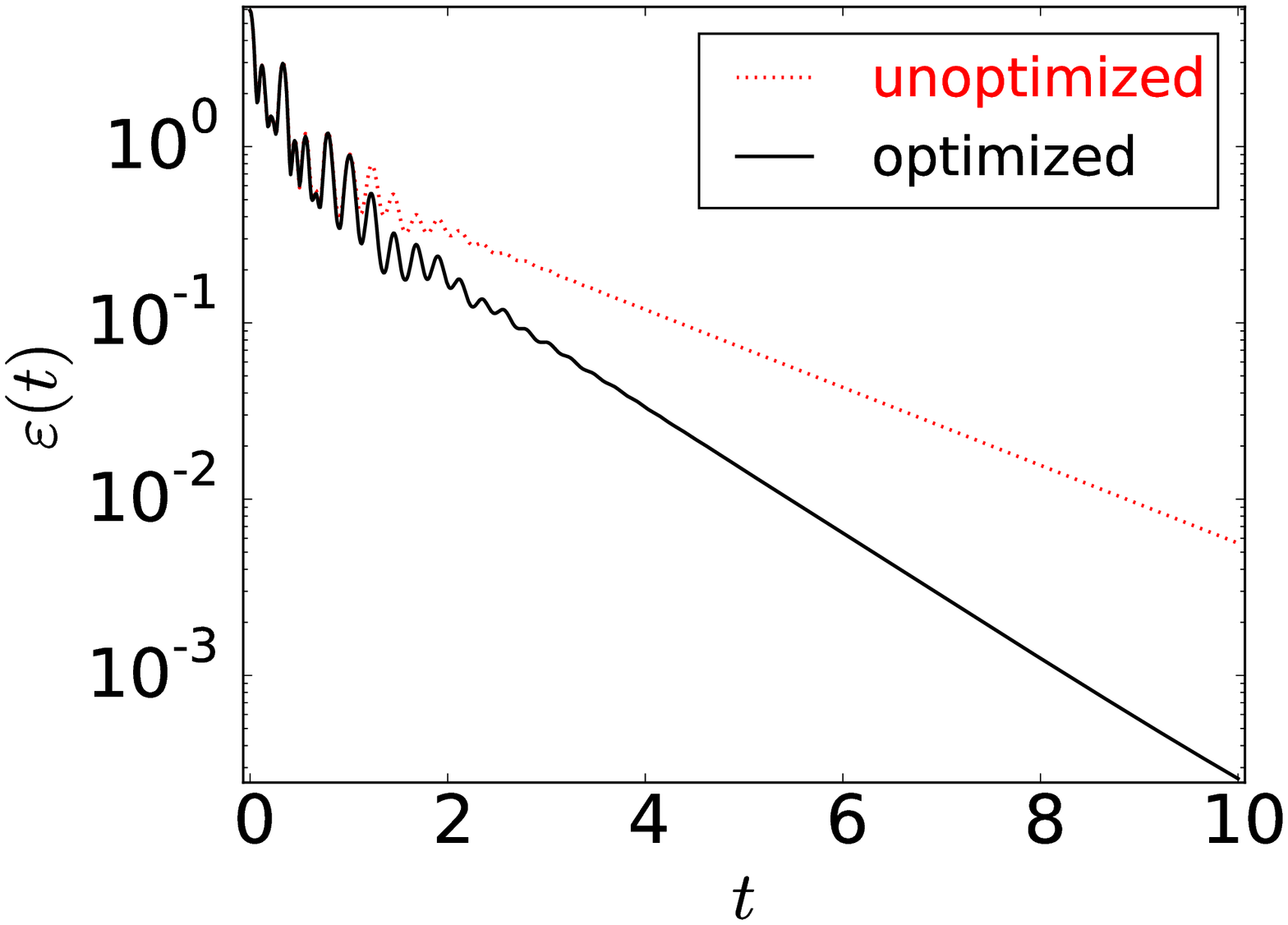}}{-1mm}{-19mm}

\caption{(Color online) (a) $\lambda_{2}$ through optimization for the RTS96 power network
under Euclidean norm constraint. (b) The unoptimized and optimized
system, where white square nodes have positive natural frequencies
(generators) while gray circular nodes have non-positive natural frequencies
(loads or relay nodes). Edge color intensity encodes $\cos(\theta_{i}^{*}-\theta_{j}^{*})$.
(c) Response of the RTS 96 power network governed by the first-order
Kuramoto model. (d) Response of the RTS 96 power network governed
by the second-order Kuramoto model with unit damping $D_{i}=1$ and small
inertia $M_{i}=0.2$. In both (c) and (d), the disturbance, drawn
from the Gaussian distribution with mean zero and standard deviation
$0.05\,\text{rad}$, was applied to the steady state of phase oscillators
at $t=0$. \label{fig:rts96_l2norm_constraint} }
\end{figure}

To obtain a non-trivial solution with optimal stability, we introduce
an additional Euclidean norm constraint,
\begin{equation}
\|\omega\|_{2}^{2}=\sum_{i}\omega_{i}^{2}\geq c,\label{eq:l2_norm_constraint}
\end{equation}
which treats all nodes in equal footing and doesn't emphasize the
role of import nodes, say, hubs. The constraint optimization is solved
by the barrier method, which is a particular interior point algorithm~\cite{boyd2004convex}.
Although the constraint Eq. (\ref{eq:l2_norm_constraint}) is nonconvex
and global optimum may not be attainable, we find in our numerical
experiments that the barrier method can efficiently achieve a satisfactory
stationary point.

In Fig.~\ref{fig:rts96_l2norm_constraint}(a) we plot the optimization
process of the RTS96 power network with constraint parameter $c=0.99\|\omega_{0}\|_{2}^{2}$,
where $\omega_{0}$ is the same as the initial natural frequency in
Fig.~\ref{fig:RTS96_unconstrained}. The corresponding unoptimized
and optimized system is shown in Fig.~\ref{fig:rts96_l2norm_constraint}(b). The edge $(318, 223)$ 
and edge $(325, 121)$ are the inter-connections between two components. In the extreme case, if both 
of them are overloaded with $|\theta_i^* - \theta_j^*|=\pi/2$ or $\cos(\theta_i^* - \theta_j^*)=0$, 
then the meta-graph with edge weight $W(\theta^*)_{ij}$ becomes disconnected into two parts, and $\lambda_2(L(\theta^*))$ 
will become zero, signaling the onset of instability of the system~\cite{ARAPOSTHATIS1981, Manik2014}. 
In our case, edge $(325, 121)$ is heavily loaded in the unoptimized system, while it is significantly 
destressed in the optimized system, achieving a more stable state as revealed by the increment of $\lambda_2(L(\theta^*))$.

To illustrate the improved stability of the optimized system related to an unoptimized one, 
we impose a small disturbance $\delta_{i}$ to the steady state at 
$t=0$, $\theta_{i}(t=0)=\theta_{i}^{*}+\delta_{i}$ and let the system
evolve according to both the first- and second-order Kuramoto model. In Figs.
\ref{fig:rts96_l2norm_constraint}(c) and ~\ref{fig:rts96_l2norm_constraint}(d)
we monitor the discrepancy between $\theta(t)$ and the steady state
$\varepsilon(t):=\sum_{i}|\theta_{i}(t)-\theta_{i}^{*}|$. It is observed
that the optimized system converges to the steady state more rapidly
than the unoptimized system.

\subsection{Properties of optimized systems}

\begin{figure}
\hspace{-2mm}
\topinset{(a)}{\includegraphics[scale=0.21]{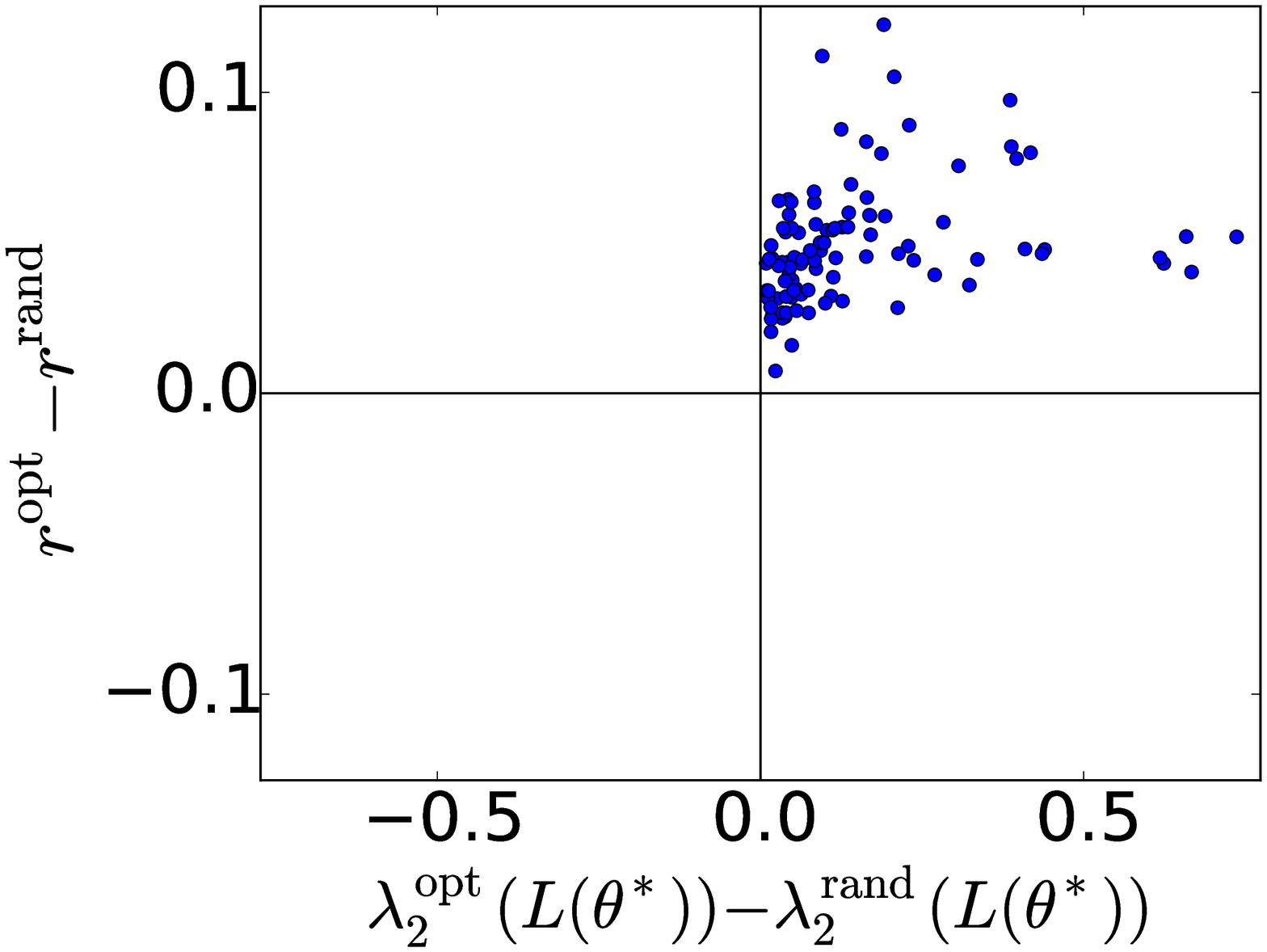}}{-1mm}{-19mm}
\topinset{(b)}{\includegraphics[scale=0.21]{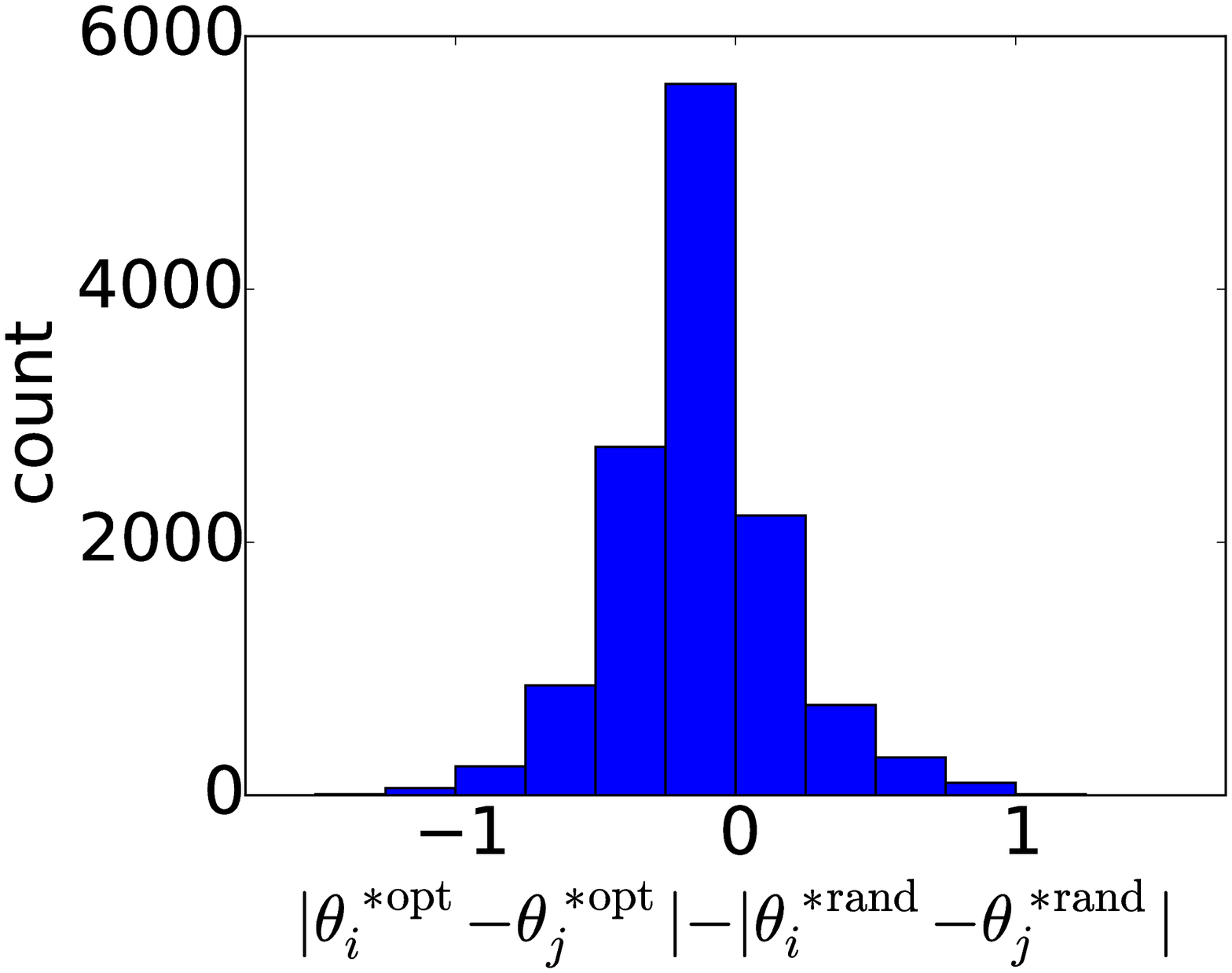}}{-1mm}{-21mm}\\
\hspace{-1mm}
\topinset{(c)}{\includegraphics[scale=0.21]{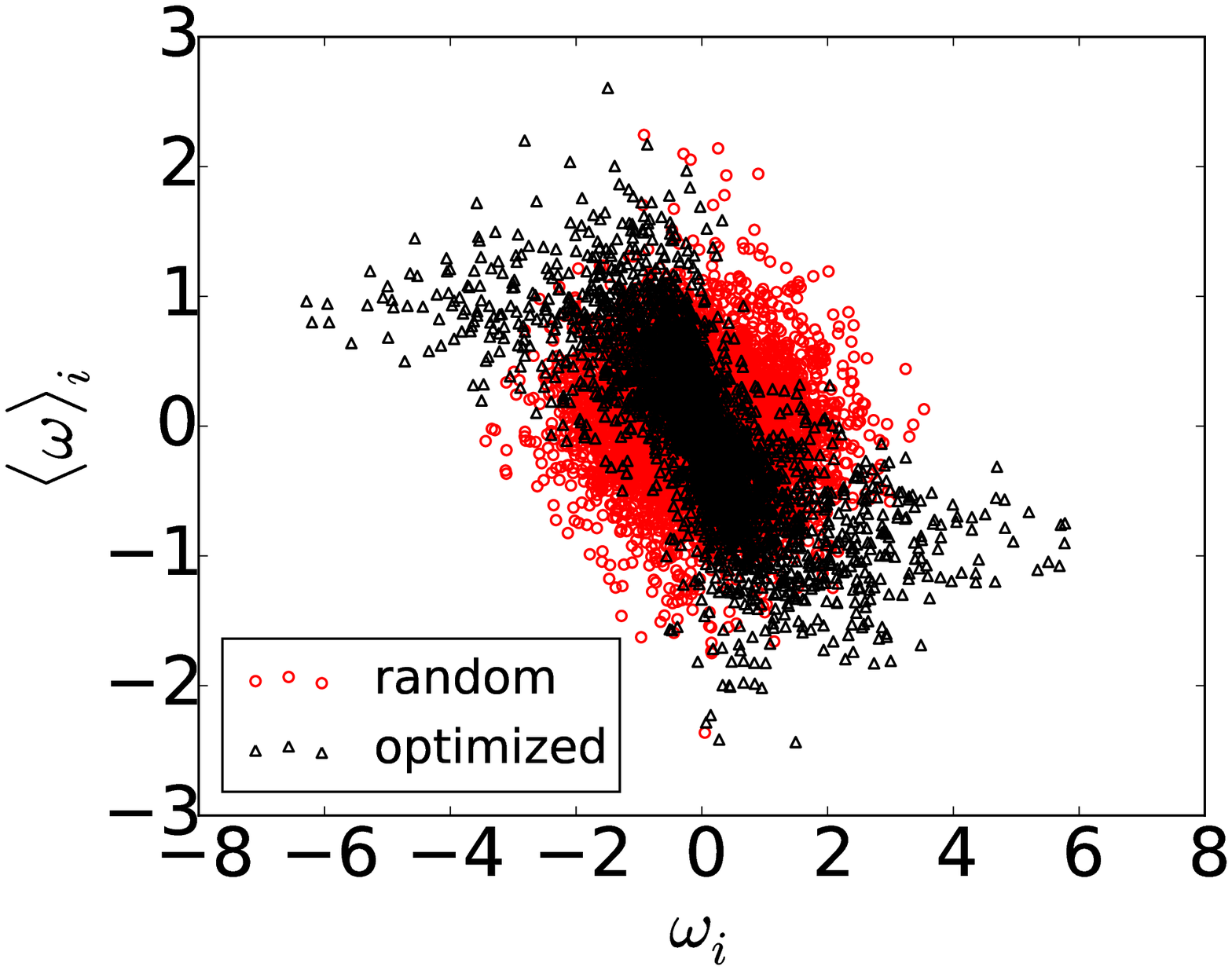}}{-1mm}{-19mm}
\topinset{(d)}{\includegraphics[scale=0.21]{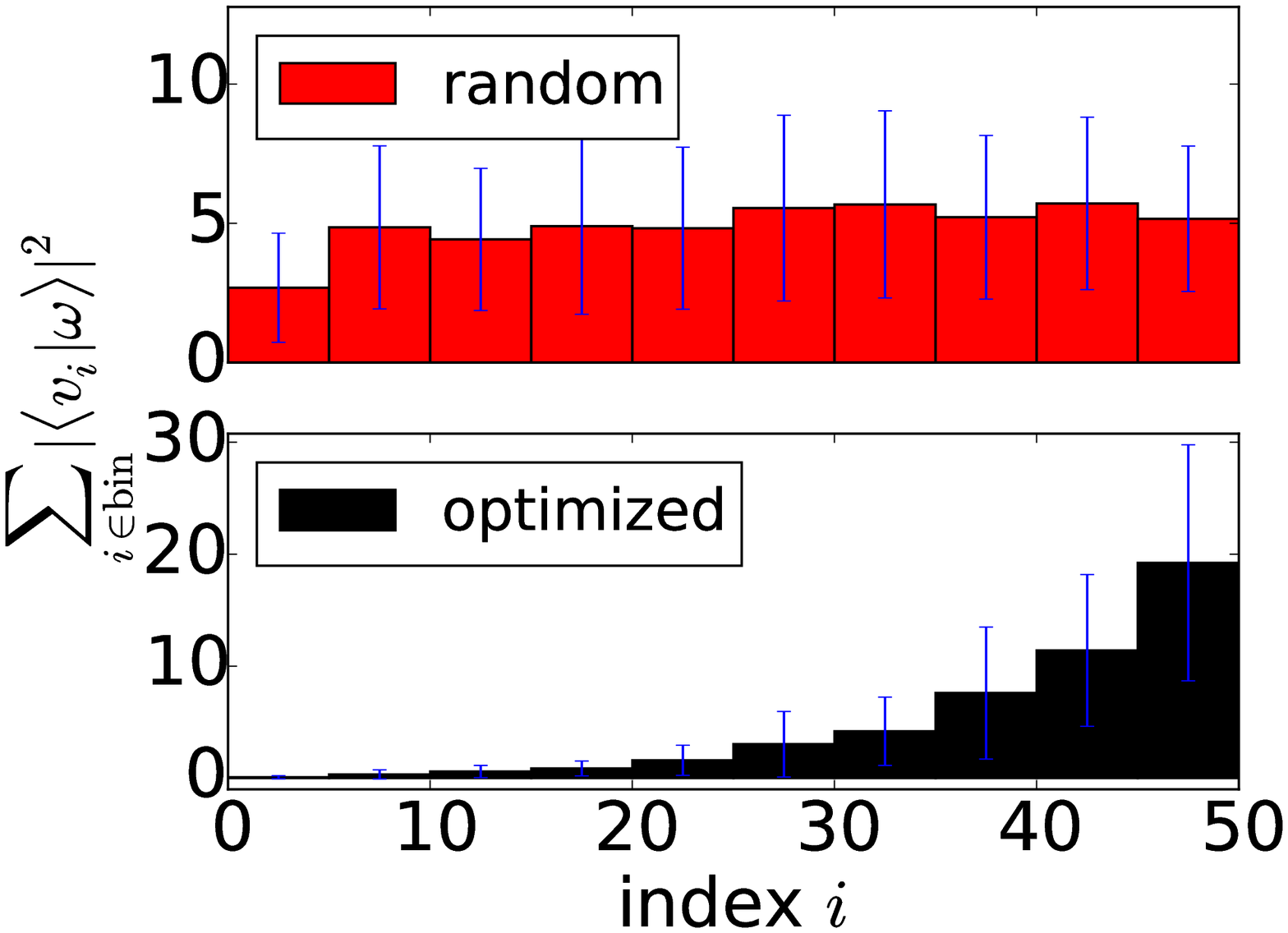}}{-1mm}{-19mm}
\caption{(Color online) Properties of the optimal system compared to systems with random frequencies.
The networks are 100 realizations of ER random graphs with 50 nodes
and edge connection probability $p=0.1$. (a) Correlation of increment of $\lambda_{2}(L(\theta^{*}))$
and order parameter $r$. (b) Histogram of changes of phase angle differences
among all edges $(i,j)$ in all realizations. (c) Average neighbor
frequency $\langle\omega\rangle_{i}=\sum_{j\in\partial i}\omega_{j}/d_{i}$
vs natural frequency $\omega_{i}$. (d) Alignments of natural frequencies
with graph Laplacian eigenvectors, i.e., $|\langle v_{i}|\omega\rangle|^{2}$
where $v_{i}$ is the normalized eigenvector corresponding to the $i$-th smallest
eigenvalue. The data was divided into 10 bins and $|\langle v_{i}|\omega\rangle|^{2}$
was first summed inside every bin for each sample, after which
the sample mean and standard deviation of the bin summation quantity
$\sum_{i\in\text{bin}}|\langle v_{i}|\omega\rangle|^{2}$ was calculated. 
\label{fig:optimal_frequency_properties}}
\end{figure}

\begin{figure*}
\topinset{(a)}{\includegraphics[scale=0.3]{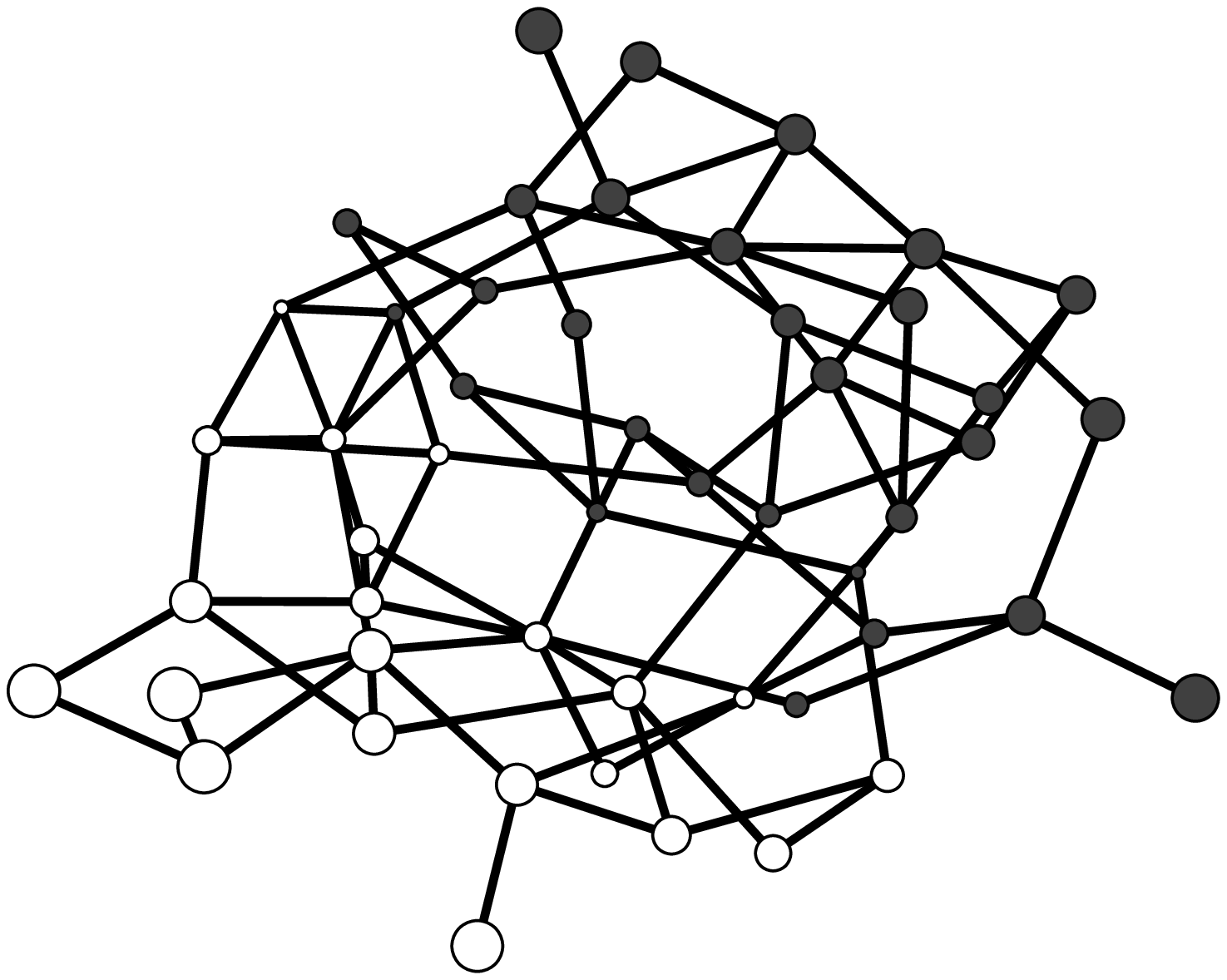}}{-1mm}{-19mm}
\topinset{(b)}{\includegraphics[scale=0.3]{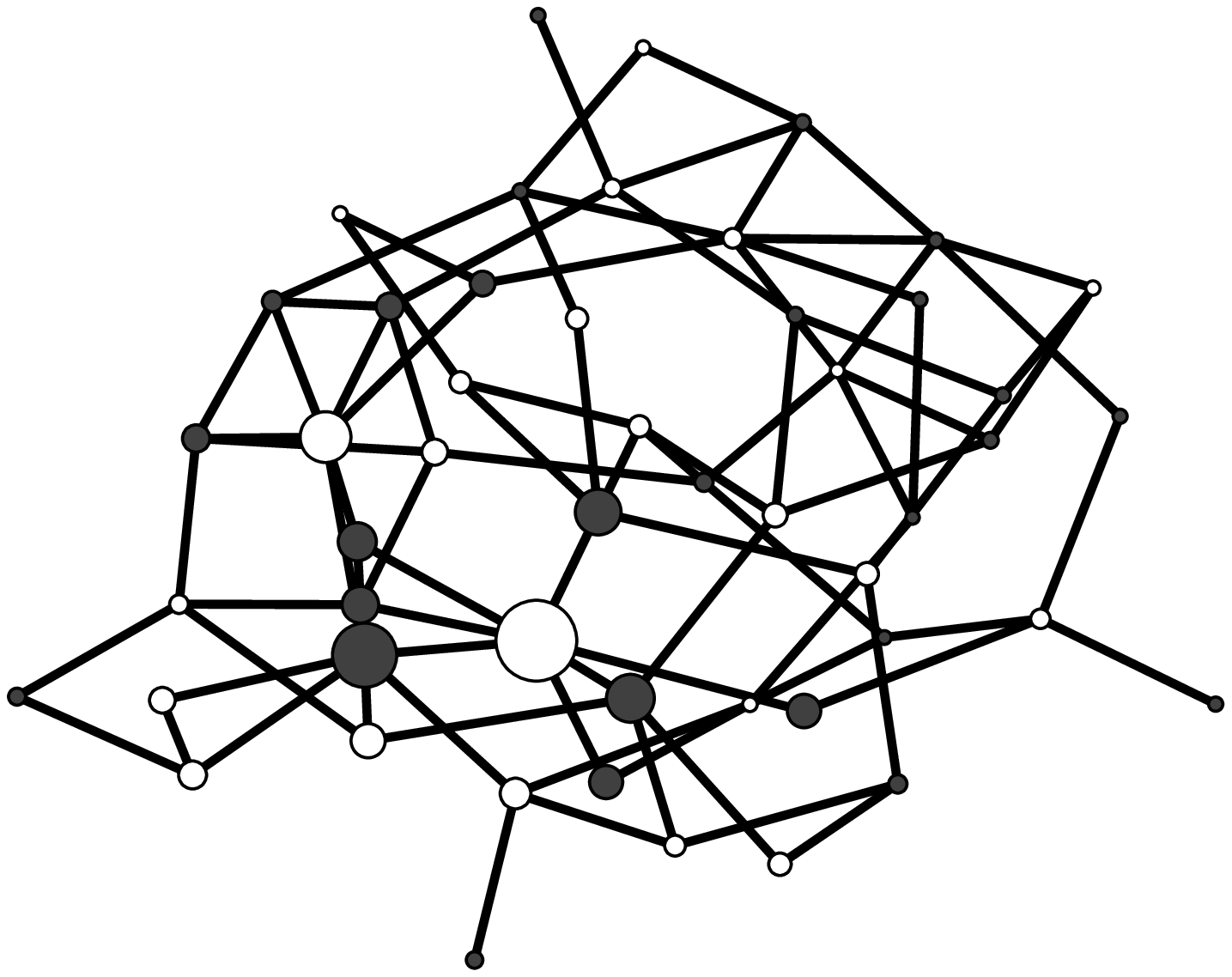}}{-1mm}{-19mm}
\topinset{(c)}{\includegraphics[scale=0.3]{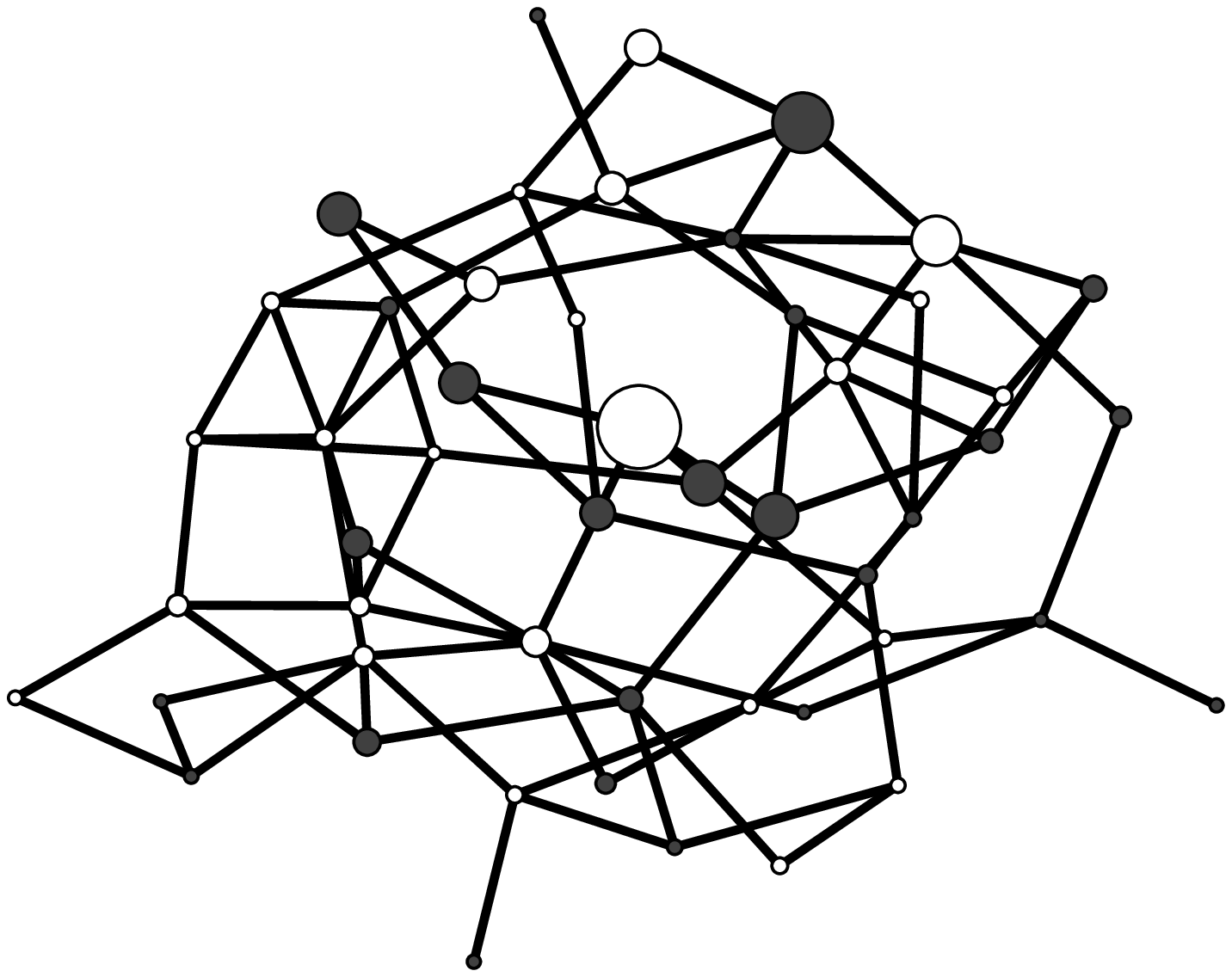}}{-1mm}{-19mm}
\caption{(a) Eigenvector $v_2$ 
corresponding to the second smallest eigenvalue $\lambda_2$ of the graph Laplacian matrix of a specific 
ER random graph, depicted on the network. The color indicates the sign of $v_{2j}$ on node $j$, i.e., 
white node corresponds to $v_{2j}>0$, while gray node corresponds to $v_{2j} \leq 0$. The size of the node 
indicates the strength of $|v_{2j}|$ on that node. (b) Eigenvector $v_N$ corresponding to the largest eigenvalue 
 $\lambda_N$ of the graph Laplacian matrix depicted on the network. (c) Frequencies $\omega^{\text{opt}}$ corresponding to the algebraic connectivity of the state-dependent Laplacian matrix evaluated at the optimal algebraic connectivity. Similarly, the white nodes correspond to $\omega^{\text{opt}}>0$, while the gray nodes correspond to $\omega^{\text{opt}} \leq 0$.}
\label{fig:ER50_eigenvectors}
\end{figure*}

In the following, we explore some general properties of the optimal
systems under the Euclidean norm constraint. The networks are ER random
graphs with 50 nodes and every pair of nodes are connected with probability
$p=0.1$. 
As found in Fig.~\ref{fig:optimal_frequency_properties}(a), not only does the optimization result in improving the objective function $\lambda_2(L(\theta^*))$, but also the Kuramoto order
parameter $r:=N^{-1}|\sum_j e^{i\theta^*_j}|$. In fact, more coherent phase angles in general imply
smaller phase angle differences $|\theta_{i}^{*}-\theta_{j}^{*}|$
and larger edge weight $W(\theta^{*})_{ij}=K_{ij}\cos(\theta_{i}^{*}-\theta_{j}^{*})$,
in which case the state-dependent network will be better connected
with a higher algebraic connectivity. Thus it is not surprising that
there is a correlation between the enhancements of $r$ and $\lambda_{2}(L(\theta^{*}))$.
We show in Fig.~\ref{fig:optimal_frequency_properties}(b) that the
decrease of phase angle differences $|\theta_{i}^{*}-\theta_{j}^{*}|$
after optimization is much more common than increase.

It is found in previous studies that natural frequencies which optimize
$r$ subject to constraint of the form $\|\omega\|_{2}^{2}=constant$
have negative correlations between neighboring frequencies, and align
with eigenvectors corresponding to large eigenvalues of graph Laplacian
\cite{Skardal2014}. We show in Figs.~\ref{fig:optimal_frequency_properties}(c) and ~\ref{fig:optimal_frequency_properties}(d)
that such properties are also observed in natural frequencies, which
optimize $\lambda_{2}(L(\theta^{*}))$. In the case of power grids on such networks, the negative correlations 
between neighboring frequencies at the optimum imply that a supply node ($\omega_i>0$) is more likely to be
connected to demand nodes ($\omega_i<0$) and vice versa. This indicates that the system stability
favors distributed power sources if all the nodes are not constrained, which is similar to the phenomenon observed 
in Ref.~\cite{Rohden2012} that decentralized power grids promote synchrony.

However, the pathways of achieving optimality with decentralized networks are different. In Ref.~\cite{Skardal2014} decentralization was achieved by maximizing the overlap of the configuration with the eigenvector of the {\em largest} eigenvalue of the graph Laplacian matrix, whereas in our work, optimal stability is achieved by maximizing the {\em smallest} positive eigenvalue of the state-dependent Laplacian matrix.

Further insight can be obtained from the alignments of optimal frequencies or power injections with the eigenvectors of \textit{graph Laplacian matrix} $L[K]$. We depict in Figs.~\ref{fig:ER50_eigenvectors}(a) and ~\ref{fig:ER50_eigenvectors}(b) the eigenvectors corresponding to the second smallest and largest eigenvalues of $L[K]$ of an ER graph, denoted as $v_2$ and $v_N$. In Fig.~\ref{fig:ER50_eigenvectors}(a), the network is partitioned into two connected subgraphs by $v_2$, with the positive components of $v_2$ belonging to one subgraph and the negative components belonging to the other, and there are only limited number of edges connecting them. It constitutes an example of graph bipartition by spectral method~\cite{Fiedler1975, Holzrichter1999}. If the power injection is aligned with $v_2$, i.e., $\omega\propto v_2$, then the implication is an extensive transportation of resources from one group to the other, as illustrated by the large phase difference across the link (325,121) in Fig.~\ref{fig:rts96_l2norm_constraint}(b), rendering the boundary between the two groups vulnerable. On the contrary, as shown in Fig. ~\ref{fig:ER50_eigenvectors}(b), the subset of positive components of $v_N$ (white) is maximally connected to the subset of negative components (gray), yielding a decentralized configuration. The observed suppression of alignment of $\omega$ with $v_2$ in Fig.~\ref{fig:optimal_frequency_properties}(d) in the optimized systems implies that the domain-wide fluctuations of resource or power is inhibited to enhance stability after optimization. On the other hand, the alignment of $\omega$ with $v_N$ is enhanced, which implies that the optimization of the system stability encourages local transmission. As shown in Fig.~\ref{fig:ER50_eigenvectors}(c), power injection on the white nodes tends to have distributed power sources.

\subsection{Difference between $\lambda_{2}(L(\theta^{*}))$ and $r$\label{sub:coupled_ER}}

\begin{figure}
\topinset{(a)}{\includegraphics[scale=0.25]{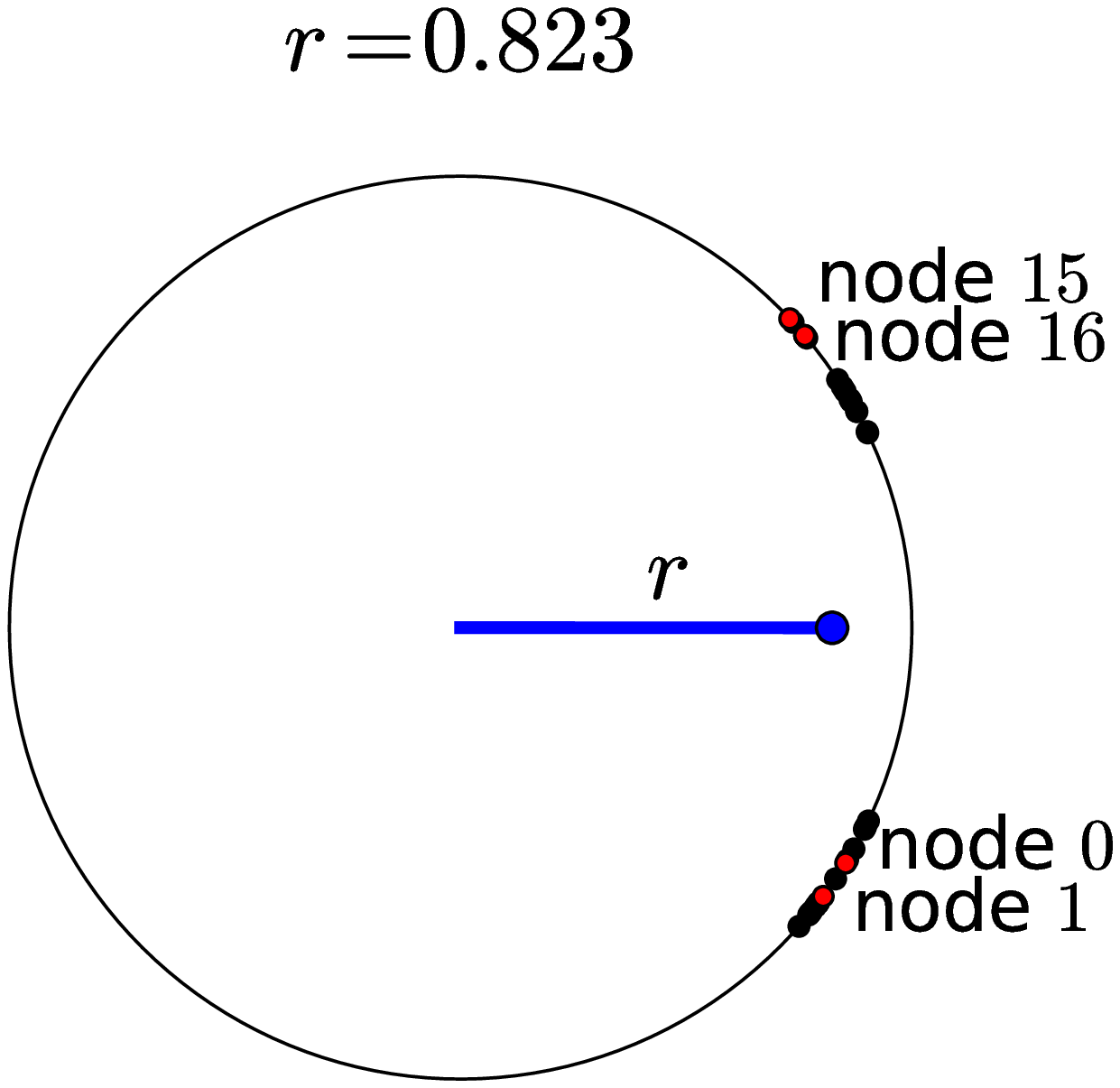}}{-1mm}{-15mm}
\topinset{(b)}{\includegraphics[scale=0.25]{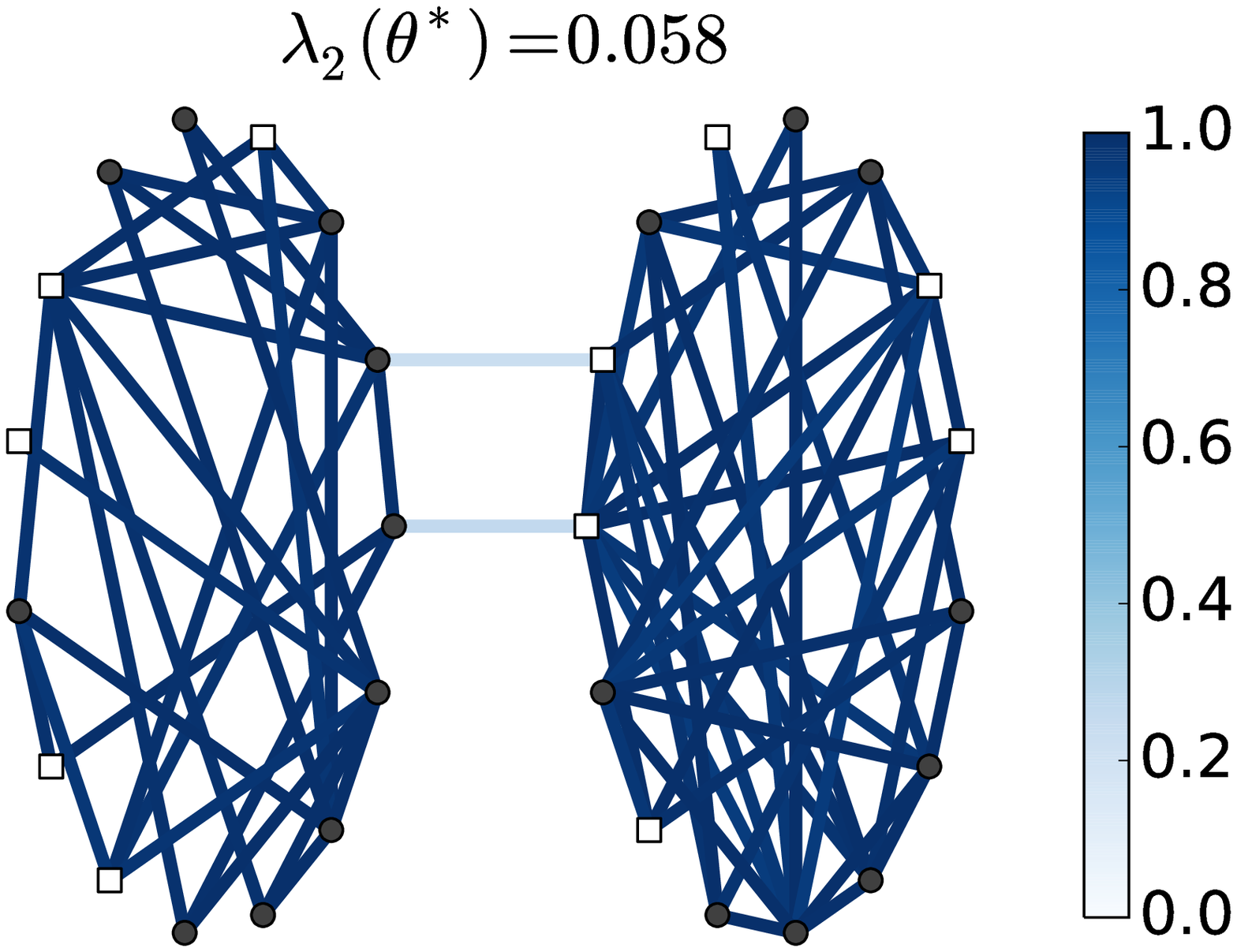}}{-1mm}{-19mm}\\
\topinset{(c)}{\includegraphics[scale=0.25]{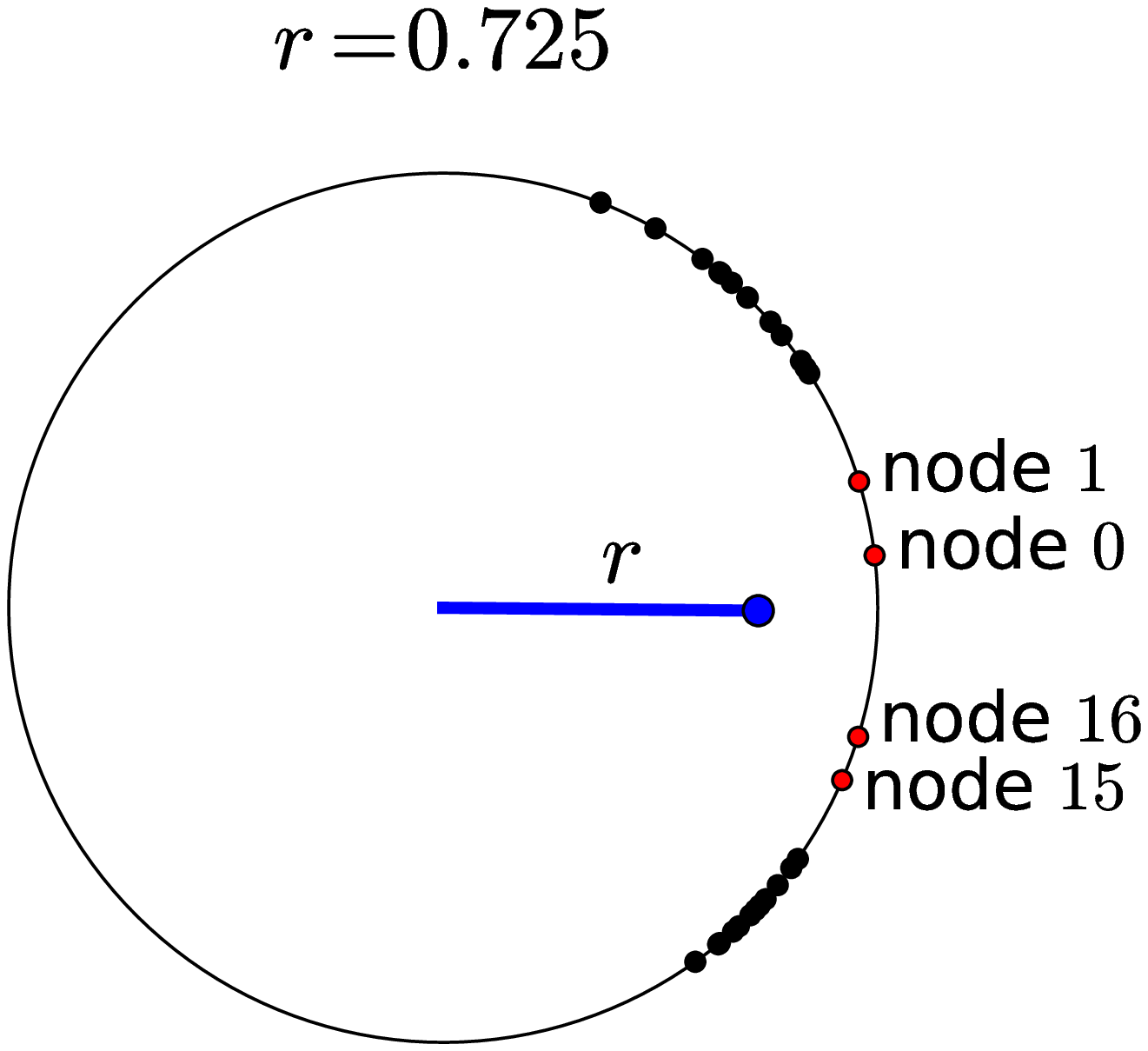}}{-1mm}{-15mm}
\topinset{(d)}{\includegraphics[scale=0.25]{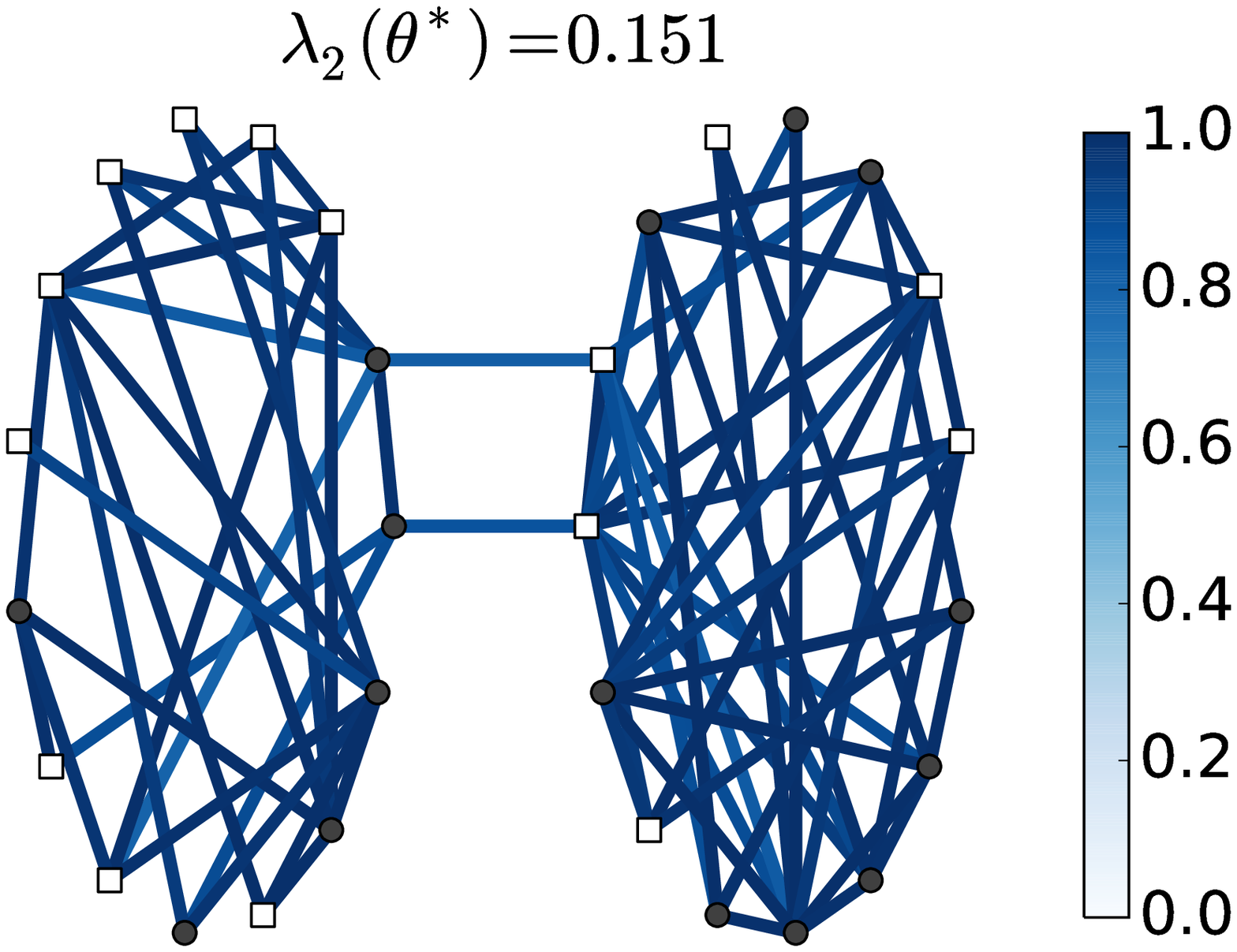}}{-1mm}{-19mm}

\caption{(Color online) Phase angles $\theta^{*}$ and state-dependent edge weights $W(\theta^{*})$
in a two-module network. In both cases, the $L^{2}$-norm of natural
frequency is $\|\omega\|_{2}=4.26$. (a) Phases of the system depicted
on the unit circle in Case 1. (b) The state-dependent edge weight
$W(\theta^{*})_{ij}=K_{ij}\cos(\theta_{i}^{*}-\theta_{j}^{*})$ in
Case 1. (c) Phases of the system depicted on the unit circle in Case
2. (d) The state-dependent edge weight $W(\theta^{*})_{ij}=K_{ij}\cos(\theta_{i}^{*}-\theta_{j}^{*})$
in Case 2. \label{fig:coupled_ER}}
\end{figure}

Observing the similarity of the results of optimizing $\lambda_{2}(L(\theta^{*}))$
with the Euclidean norm constraint and those of optimizing $r$ with the same constraint, it is
tempting to conclude that the more synchronized a system the more stable it is and
one can improve the system stability by just increasing the order
parameter $r$, which can be much simpler. However, we argue that
while such a judgment is valid in many cases like the above homogeneous ER
graphs, it is not necessarily a universal rule. In most cases, optimizing $r$ will not be the most
efficient way to enhance the system stability. Moreover, there is a conceptual difference
between the two quantities. The Kuramoto order parameter $r$ is a
measure of coherence of phase angles of all oscillators in a global
and average sense, which cannot identify the role of critical edges
in maintaining stability, e.g., the interconnections between modules.
To be more concrete, we consider a simple network which is composed
of two modules, each corresponding to a small random graph, as sketched
in Fig.~\ref{fig:coupled_ER}(b). The coupling of each edge is set
to be $K_{ij}=1$.

In Case 1, we suppress the intra-module transportation and encourage
the inter-module transportation, which leads to phases that are coherent
inside each module but have a large separation between the two modules,
as shown in Fig.~\ref{fig:coupled_ER}(a). The phase coherence	
inside each module leads to a relatively high Kuramoto order parameter
$r=0.823$. However, the large inter-module phase difference indicates
the edge $(0,15)$ and edge $(1,16)$ are highly stressed with a low
state dependent edge weight $W(\theta^{*})_{ij}=K_{ij}\cos(\theta_{i}^{*}-\theta_{j}^{*})$,
resulting in a small state algebraic connectivity $\lambda_{2}(L(\theta^{*}))=0.058$
as shown in Fig.~\ref{fig:coupled_ER}(b). In Case 2, the system is
perturbed and the phases become more dispersed, leading to a smaller
Kuramoto order parameter $r=0.725$. But the phase differences along
edge $(0,15)$ and edge $(1,16)$ are much reduced. This significantly
increases the edge weights $W(\theta^{*})_{ij}$ of these two edges and hence the state algebraic
connectivity reaches $\lambda_{2}(L(\theta^{*}))=0.151$, since edge
$(1,15)$ and edge $(1,16)$ are the inter-module connections whose
edge weights are crucial for the algebraic connectivity. This simple
example highlights the essence of using $\lambda_{2}(L(\theta^{*}))$
as a cost function for measuring stability in general networks.

\subsection{Inclusion of practical power grid constraints}

The Euclidean norm-constrained optimization problem above treats all
nodes on equal footing where a supplier can become a consumer and
vice versa. This will not be realistic if we consider power grid
applications. In this section, we consider two problems regarding
practical constraints of power grid operations. 

In Problem 1, both the supply and the demand are restricted to vary within
a certain range. Furthermore, regulating both the generation and consumption may be necessary in future grids with the introduction of renewable energy. Hence specifically we consider the constraint $\omega_{0i}-\alpha|\omega_{0i}|\leq\omega_{i}\leq\omega_{0i}+\alpha|\omega_{0i}|$
for $i$ to be either a supply node or demand node, where $\omega_{0i}$
is the natural frequency of the original system and the parameter
$\alpha$ satisfies $0<\alpha\leq1$. For the relay node with $\omega_{0i}=0$,
the natural frequency will remain unchanged throughout optimization
$\omega_{i}=\omega_{0i}=0$. 

In Problem 2, only the supply nodes
with $\omega_{0i}>0$ are allowed to schedule their productions with fraction
$\alpha$, while the demands must be satisfied and the relay nodes should
also be fixed, i.e., $\omega_{i}=\omega_{0i}$ for $\omega_{0i}\leq0$.
To deal with both the inequality and equality constraints, the primal-dual
interior point method in convex optimization is applied in these problems.
Although we always make the supply and demand balanced in every iteration,
we discovered that imposing the additional constraint $\sum_{i}\omega_{i}=0$
into the definition of the problem can significantly facilitate the
convergence of the algorithm.

\begin{figure}
\hspace{-2mm}
\topinset{(a)}{\includegraphics[scale=0.22]{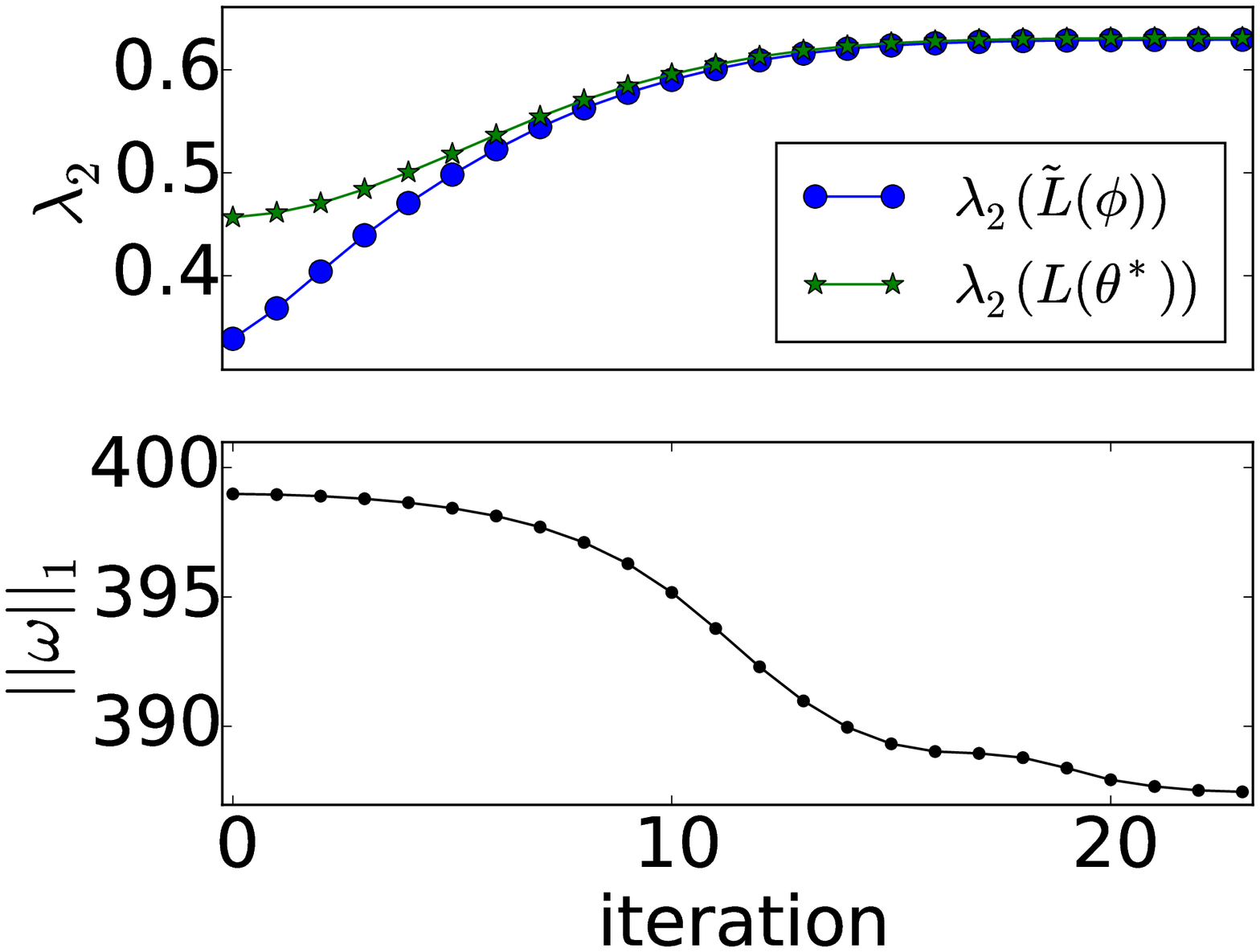}}{1mm}{-10mm}
\topinset{(b)}{\includegraphics[scale=0.22]{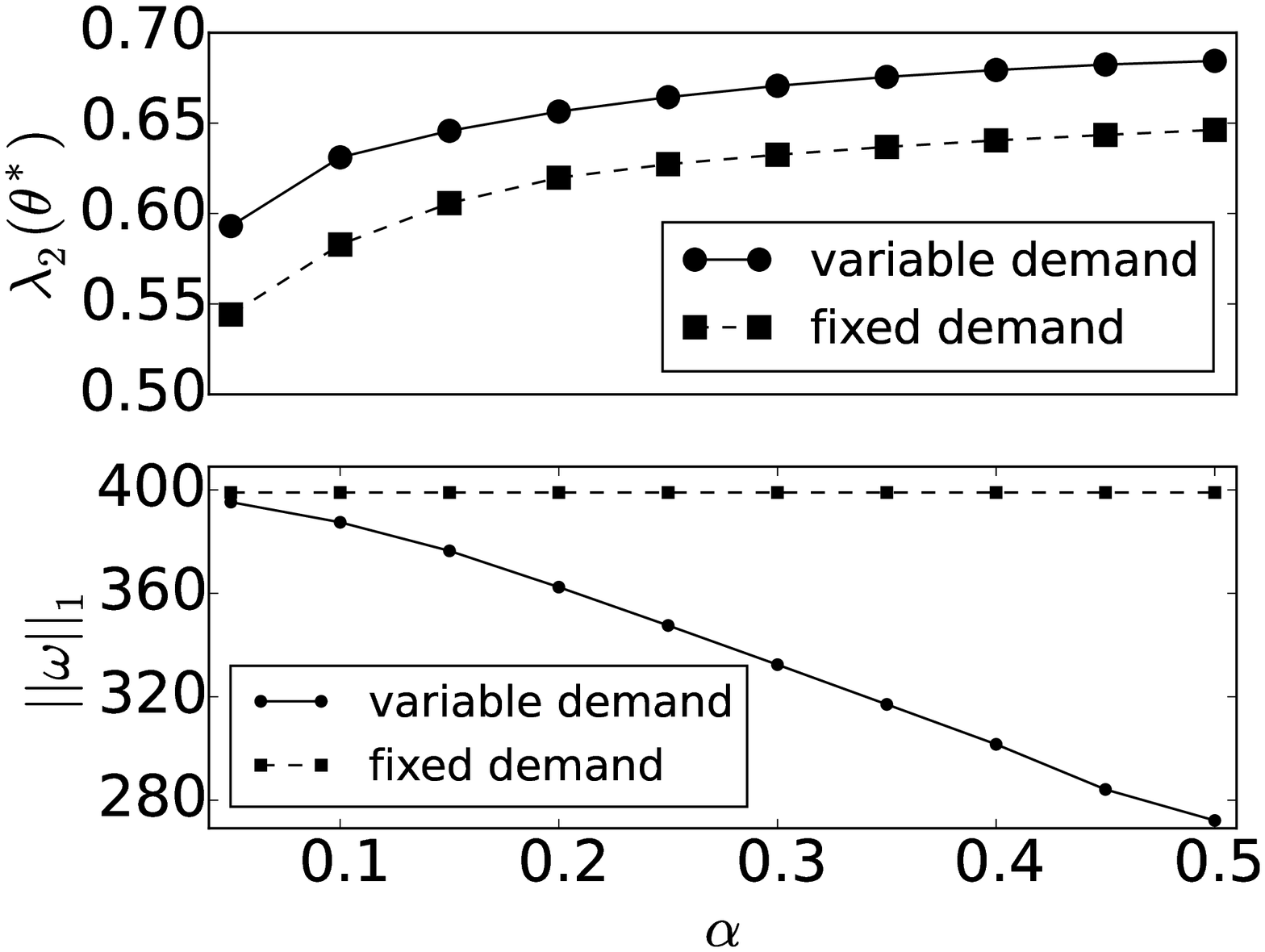}}{2mm}{-11mm}

\caption{(Color online) (a) $\lambda_{2}$ and $\|\omega\|_{1}$ through optimization for
the RTS96 power network under linear constraints at $\alpha=0.1$. The $L_{1}$-norm
of natural frequency $\|\omega\|_{1}:=\sum_{i}|\omega_{i}|$ is twice
of total generation or total consumption. (b) $\lambda_{2}(L(\theta^{*}))$
and $\|\omega\|_{1}$ of the optimal system as a function of $\alpha$
with variable demand (Problem 1) and fixed demand (Problem 2).\label{fig:RTS96_power_grid_constraint}}
\end{figure}

In Fig.~\ref{fig:RTS96_power_grid_constraint}(a), we plot the optimization
process of the RTS96 power network with constraints of Problem 1.
The primal-dual interior point algorithm can bring the system to optimum
effectively. We also monitor the $L^1$-norm of $\omega$, defined as $\|\omega\|_1:=\sum_i|\omega_i|$, 
which is twice the total production or total consumption. During optimization, the system is also destressed as
indicated by the decrement of $\|\omega\|_{1}$. 
In Fig.~\ref{fig:RTS96_power_grid_constraint}(b),
we plot $\lambda_{2}(L(\theta^{*}))$ and $\|\omega\|_{1}$ as a function
of $\alpha$ with constraints of both Problem 1 and Problem 2. It
is observed that $\lambda_{2}(L(\theta^{*}))$ increases with $\alpha$
for both cases with variable demands and fixed demands. This is not
surprising since the feasible region of the problem with larger $\alpha$
is a superset of the one with smaller $\alpha$, and a larger feasible
region gives the system more flexibility to search for more stable
state. The system can achieve higher stability with variable demands
in Problem 1 than the fixed demand in Problem 2, which is also due
to more degrees of freedom to vary in Problem 1. Our method can solve
both problems satisfactorily. 

\subsection{Behavior at optimal coupling strengths}

\begin{figure}
\topinset{(a)}{\includegraphics[scale=0.22]{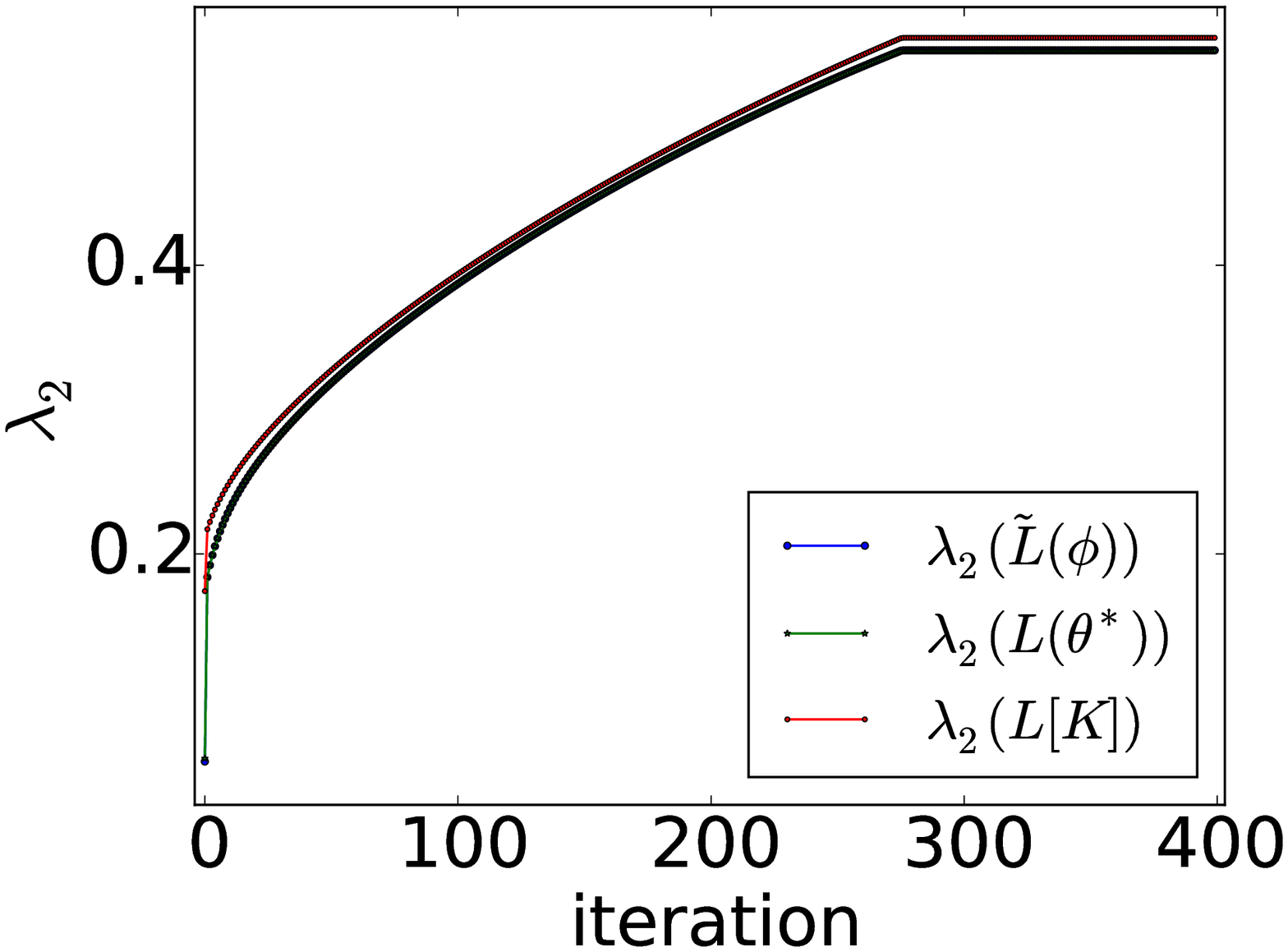}}{-1mm}{-18mm}
\topinset{(b)}{\includegraphics[scale=0.24]{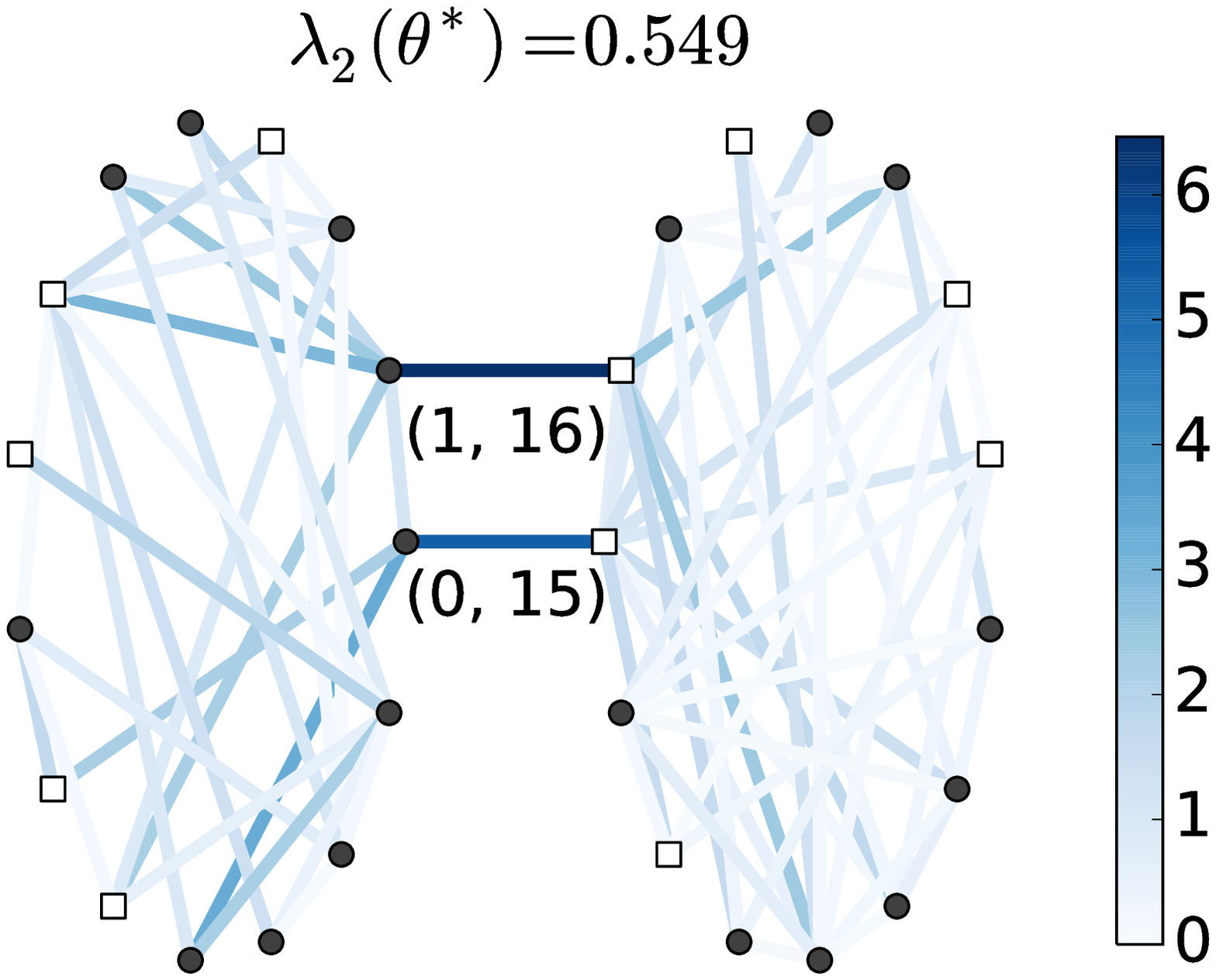}}{-1mm}{-18mm}

\caption{(Color online) Optimizing the state algebraic connectivity by updating coupling strengths.
(a) State algebraic connectivity $\lambda_{2}(L(\theta^{*}))$, $\lambda_{2}(\tilde{L}(\phi))$,
and graph algebraic connectivity $\lambda_{2}(L[K])$ through optimization.
The initial state is the same as Case 1 in Sec.~\ref{sub:coupled_ER}.
(b) The state-dependent edge weight $W(\theta^{*})_{ij}=K_{ij}\cos(\theta_{i}^{*}-\theta_{j}^{*})$
in the optimal state. Note the scale of color code is different from
the cases of Fig.~\ref{fig:coupled_ER}.\label{fig:Update_K} }
\end{figure}

Lastly, we consider behavior at the optimal state algebraic connectivity by
updating the coupling strengths. To avoid indefinite solutions, we
impose a simple constraint
\begin{equation}
\sum_{(i,j)}K_{ij}=K_{\text{total}},\label{eq:total_K_constraint}
\end{equation}
where $K_{\text{total}}$ represents the availability of the total capacity,
and $K_{ij}$ is constrained to be non-negative. Due to the high complexity
of computing the Hessian, we only consider the gradient ascent update.
To preserve the resource constraint, the approximated gradient $\nabla_{K}\lambda_{2}(\tilde{L}(\phi))$
as calculated by Eq. (\ref{eq:gradient_wrt_K}) is projected onto
the feasible region, after which the coupling strengths are updated.
In Fig.~\ref{fig:Update_K}(a), we plot the optimization process of
the projected gradient update on the two-module network discussed
in Sec.~\ref{sub:coupled_ER}, and the initial condition is the
same as Case 1 in Sec.~\ref{sub:coupled_ER}. It is shown that
redistributing the coupling strengths can significantly improve both
the graph algebraic connectivity and state-algebraic connectivity,
reaching a more stable state. In Fig.~\ref{fig:Update_K}(b), we sketch
the state-dependent edge weight in the optimal state. Contrary
to the un-optimized system in Fig.~\ref{fig:coupled_ER}(b), the optimized
system exhibits large edge weight $W(\theta^{*})_{ij}$ in edge $(1,16)$
and edge $(0,15)$, the interconnections between the two modules,
which favors higher state algebraic connectivity. For each module,
the nodes are well connected and the need for transporting resource
is modest. Thus, the coupling strengths inside each module are sacrificed
so that the system can invest more on the the critical edges. The effects of increasing stability 
by investments on the inter-area links are also studied and demonstrated in Ref.~\cite{Wang2016}, where the interlinks
are added one by one according to the greedy search strategy instead of updating the existing links as in
our approach. These phenomena highlight the importance of strengthening the inter-connections between
different communities of the grid.  

\section{Discussion}

In this paper, we studied the optimization of synchronization stability
of the Kuramoto model by updating the natural frequencies or coupling
strengths. The proposed cut-set space approximation can accurately
estimate the network flows of steady states and thus simplify the
objective function, i.e., the state algebraic connectivity whose increment
can increase the stability of the phase-locked steady states of both the first- and second-order Kuramoto model. 
Such an approximation leads to compact expressions of gradient and Hessian of the cost function.
Together with the interior point algorithm or projected gradient ascent,
our method can cope with various constraints, which is shown to be
effective and efficient. There is a general correlation between the
optimization of the Kuramoto order parameter and the state algebraic
connectivity, especially in the homogeneous networks. However, the
Kuramoto order parameter cannot represent the role of critical links,
e.g., inter-module connections, which is crucial to the synchronization
stability. In light of this consideration, the state algebraic connectivity
is a more appropriate cost function for the measure of stability.
Our framework has potential applications in improving the stability
of power grids which are usually simplified to a second-order Kuramoto
model. The method also sheds light on the treatments of general nonlinear 
eigenvalue optimization problems.

Nevertheless, there are many other aspects to consider concerning
the application of power grids, such as extending our formalism to
non-uniform inertia or damping, lossy transmissions, effect of changes
of network topology due to breakdown of grid elements, etc. In addition,
our method is based on the assumption of non-degenerate state algebraic
connectivity, which may not hold in highly symmetric networks, and
how to achieve an optimum under general constraints in these networks
remains to be explored. Last, our study considers only linear stability
which assumes small disturbances. While we found that the decentralized configuration has optimal stability for small disturbances, there were indications that decentralization may reduce the dynamic stability for moderate perturbations~\cite{Rohden2012}. This may require us to adopt an augmented objective function in future studies. The recently developed basin stability
approach~\cite{Menck2014} can be complementary to our approach, and the combination
of the two views may be able to provide more comprehensive understanding
of the system stability.

\section*{Acknowledgments}

We are grateful to D. Saad, H. Wang, P. Choi, M. Yan, H. Tsang and X. Huang for fruitful discussions. This work is supported by grants from the Research Grants Council of Hong Kong (Grants No. 605813 No. 16322616).

\bibliographystyle{unsrt}
\bibliography{ref}

\end{document}